\documentclass[10pt,oneside,english]{amsart}
\usepackage[T1]{fontenc}
\usepackage[latin9]{inputenc}
\usepackage{amssymb}

\makeatletter
\numberwithin{equation}{section} 
\numberwithin{figure}{section} 
  \@ifundefined{theoremstyle}{\usepackage{amsthm}}{}
  \theoremstyle{plain}
  \newtheorem{thm}{Theorem}[section]
  \theoremstyle{definition}
  \newtheorem{defn}[thm]{Definition}
  \theoremstyle{plain}
  \newtheorem{lem}[thm]{Lemma}
  \theoremstyle{plain}
  \newtheorem{cor}[thm]{Corollary}
  \theoremstyle{remark}
  \newtheorem{rem}[thm]{Remark}
 \theoremstyle{definition}
  \newtheorem{example}[thm]{Example}
  \theoremstyle{remark}
  \newtheorem*{acknowledgement*}{Acknowledgement}

\makeatother

\usepackage{babel}

\begin{document}

\title{Wandering vectors and\\
the reflexivity of free semigroup algebras}

\author{Matthew Kennedy}

\address{Pure Math. Dept. \\
 University of Waterloo \\
 Waterloo, Ontario, Canada N2L 3G1}

\email{m3kennedy@uwaterloo.ca}

\begin{abstract}
A free semigroup algebra $\mathcal{S}$ is the weak-operator-closed
(non-self-adjoint) operator algebra generated by $n$ isometries with
pairwise orthogonal ranges. A unit vector $x$ is said to be wandering
for $\mathcal{S}$ if the set of images of $x$ under non-commuting
words in the generators of $\mathcal{S}$ is orthonormal.

We establish the following dichotomy: either a free semigroup algebra
has a wandering vector, or it is a von Neumann algebra. Consequences
include that every free semigroup algebra is refl{}exive, and that
certain free semigroup algebras are hyper-reflexive with a very small
hyper-reflexivity constant.
\end{abstract}

\subjclass[2000]{Primary 47L80; Secondary 47L55, 47A15}

\thanks{Research supported by NSERC}

\maketitle

\section{Introduction}

A \emph{free semigroup algebra} $\mathcal{S}$ is the weak-operator-closed
(non-self-adjoint) algebra generated by $n$ isometries $S_{1},...,S_{n}$
on a Hilbert space $\mathcal{H}$ which have pairwise orthogonal ranges,
or equivalently, which satisfy \[
S_{i}^{*}S_{j}=\begin{cases}
I & \quad\mbox{if }i=j,\\
0 & \quad\mbox{otherwise}.\end{cases}\]
Although $n$ can be finite or infinite, for notational convenience
we treat $n$ as finite and make note of any issues that arise. We
say that $S=[S_{1}\enspace...\enspace S_{n}]$ is a \emph{row isometry},
since $S$ is isometric as a row operator from $\mathcal{H}^{n}$
to $\mathcal{H}$.

Row isometries arise throughout operator theory. A theorem of Frazho,
Bunce, and Popescu shows that $n$ operators $A_{1},...,A_{n}$ which
satisfy $\sum A_{k}A_{k}^{*}\leq I$ can be \emph{dilated} to a row
isometry $S=[S_{1}\enspace...\enspace S_{n}]$ such that\[
S_{k}=\left(\begin{array}{cc}
A_{k} & 0\\
* & *\end{array}\right).\]
This is a non-commutative multivariable analogue of the Sz.-Nagy dilation
theorem.

Popescu \cite{MR1343719} showed that the norm-closed algebra generated
by any row isometry of size $n$ is completely isometrically isomorphic
to the \emph{non-commutative disk algebra} $\mathcal{A}_{n}$, and
it is well known that the C{*}-algebra generated by a row isometry
of size $n$ is isomorphic to the Cuntz algebra $\mathcal{O}_{n}$
if $\sum A_{k}A_{k}^{*}=I$, and otherwise is isomorphic to the Cuntz-Toeplitz
algebra $\mathcal{E}_{n}$. By contrast, the weak-operator-closed
algebras generated by distinct row isometries can be dramatically
different (see for example \cite{MR1823866}).

In some sense then, it is natural to study row isometries by looking
at the free semigroup algebras they generate. This idea, and with
it the definition of a free semigroup algebra, was introduced by Davidson
and Pitts \cite{MR1665248}. They observed that free semigroup algebras
often contain interesting information about the unitary invariants
of their generators.

The prototypical example of a free semigroup algebra is the \emph{non-commutative
analytic Toeplitz algebra} generated by the left regular representation
of the free semigroup on $n$ letters. This algebra, which we denote
by $\mathcal{L}_{n}$, was first studied by Popescu \cite{MR1129595}
in the context of non-commutative multivariable dilation theory. 

For $n=1$, $\mathcal{L}_{n}$ is the familiar algebra of analytic
Toeplitz operators, which is singly generated by the unilateral shift.
For $n\geq2$, $\mathcal{L}_{n}$ is no longer commutative, but it
turns out that a number of classical results about the analytic Toeplitz
operators have straightforward generalizations to this setting. This
is a large part of the motivation for the name ``non-commutative analytic
Toeplitz algebra.''

The role of $\mathcal{L}_{n}$ is of central importance in the general
theory of free semigroup algebras, and it turns out to be desirable
to isolate ``$\mathcal{L}_{n}$-like'' behavior. A free semigroup
algebra is said to be of \emph{type L} if it is algebraically isomorphic
to $\mathcal{L}_{n}$. It is important to emphasize the word ``algebraically''
here. Examples have been constructed (see for example \cite{MR1823866})
of free semigroup algebras which are of type L, and so behave algebraically
like $\mathcal{L}_{n}$, but which have a very different spatial structure.

The general structure theorem for free semigroup algebras \cite{MR1823866}
shows that every free semigroup algebra can be decomposed into $2\times2$
block-lower-triangular form, where the left column is a slice of a
von Neumann algebra, and the bottom-right entry is a type L free semigroup
algebra. It is well known (see for example \cite{MR0048700}) that
the weak-operator-closed algebra generated by a single isometry can
be self-adjoint. Davidson, Katsoulis, and Pitts \cite{MR1823866}
asked whether it was possible for a free semigroup algebra on $2$
or more generators to be self-adjoint, and some time later Read \cite{MR2186356}
(see also \cite{MR2204288}) answered in the affirmative by showing
that $\mathcal{B}(\mathcal{H})$ was a free semigroup algebra.

A notion of fundamental importance is that of a wandering vector.
A unit vector $x$ is said to be \emph{wandering} for a free semigroup
algebra $\mathcal{S}$ if the set of images of $x$ under non-commuting
words in the generators of $\mathcal{S}$ forms an orthonormal set.
It is known (see for example \cite{MR1665248}) that the spatial structure
of $\mathcal{L}_{n}$ is completely determined by the existence of
a large number of wandering vectors.

It is easy to see that the restriction of any free semigroup algebra
to the cyclic subspace generated by a wandering vector is unitarily
equivalent to $\mathcal{L}_{n}$, and so in particular is of type
L. It has been an open question for some time, however, whether every
type L free semigroup algebra necessarily has a wandering vector.
It turns out that this question is equivalent to the question of whether
every free semigroup algebra is reflexive. This can be shown using
the general structure theorem for free semigroup algebras: since every
von Neumann algebra is reflexive, the reflexivity of a free semigroup
algebra depends on the reflexivity of its type L part.

The purpose of this paper is to prove that every type L free semigroup
algebra has wandering vectors, and hence to prove that every free
semigroup algebra is reflexive.

Our approach is very much in the spirit of the ``dual algebra arguments''
which have been used with great success by Bercovici, Foias, Pearcy
and many others (see for example \cite{MR787041}), and which are
based on Brown's proof of the existence of invariant subspaces for
subnormal operators \cite{MR511974}. The fundamental idea at the
heart of these arguments is that it is often possible to prove the
existence of invariant subspaces for a weak{*}-closed  operator algebra
by showing that, in an appropriate sense, the predual of the algebra
is small.

Typically, these arguments are employed in a commutative setting,
where certain spectral and function-theoretic tools are available.
In the present non-commutative context, we rely instead on various
operator-theoretic techniques.

A clue that it might be possible to attack the present problem using
dual algebra techniques came from a recent paper of Bercovici \cite{MR1641578},
who used them to establish the hyper-reflexivity of a class of algebras
which includes the non-commutative analytic Toeplitz algebra on two
or more generators. The hyper-reflexivity of this algebra had already
been shown by Davidson and Pitts \cite{MR1665248}, with an upper
bound of $51$ on the hyper-reflexivity constant, but Bercovici's
approach yielded a surprisingly low upper bound of $3$.

Motivated by Bercovici's result, once we have shown that every type
L free semigroup algebra has a wandering vector, we go further and
show that every type L free semigroup algebra on two or more generators
is hyper-reflexive with hyper-reflexivity constant at most $3$.

\section{Preliminaries}

Let $\mathbb{F}_{n}^{+}$ denote the free semigroup in $n$ non-commuting
letters $\{1,...,n\}$, including the empty word $\varnothing$. For
a word $w$ in $\mathbb{F}_{n}^{+}$, let $|w|$ denote its length,
and let $\mathbb{F}_{n}^{k}$ denote the set of all words in $\mathbb{F}_{n}^{+}$
of length at most $k$.

Let $\mathcal{F}_{n}$ denote the ``Fock'' space $\mathcal{F}_{n}=\ell^{2}(\mathbb{F}_{n}^{+})$
with orthonormal basis $\{\xi_{w}:w\in\mathbb{F}_{n}^{+}\}$ consisting
of words in $\mathbb{F}_{n}^{+}$. For each $v$ in $\mathbb{F}_{n}^{+}$,
define an isometry $L_{v}$ by\[
L_{v}\xi_{w}=\xi_{vw},\quad w\in\mathbb{F}_{n}^{+}.\]
The map $v\to L_{v}$ gives a representation of $\mathbb{F}_{n}^{+}$,
called the \emph{left regular representation}.

The isometries $L_{1},...,L_{n}$ have pairwise orthogonal ranges.
The free semigroup algebra they generate, denoted by $\mathcal{L}_{n}$,
is called the \emph{non-commutative analytic Toeplitz algebra}. For
$n=1$, $\mathcal{L}_{n}$ is the classical analytic Toeplitz algebra,
but for $n\geq2$, $\mathcal{L}_{n}$ is no longer commutative.

We require a result for $\mathcal{L}_{n}$ which generalizes a classical
result about the analytic Toeplitz operators. An element in $\mathcal{L}_{n}$
is said to be \emph{inner} if it is an isometry, and \emph{outer}
if it has dense range. It was shown in \cite{MR1665248} that an arbitrary
element $A$ in $\mathcal{L}_{n}$ can be written as $A=BC$, where
$B$ is inner and $C$ is outer. This generalizes the classical inner-outer
factorization for elements in the analytic Toeplitz algebra.

Every element $A$ in $\mathcal{L}_{n}$ is completely determined
by its Fourier series \[
A\sim\sum_{w\in\mathbb{F}_{n}^{+}}a_{w}L_{w},\]
which is a formal power series with coefficients in $\mathcal{L}_{n}$,
where \[
A\xi_{\varnothing}=\sum_{w\in\mathbb{F}_{n}^{+}}a_{w}\xi_{w}.\]
For $k\geq1$, define the $k$-th Cesaro sum of the Fourier series
of $A$ by\[
\Gamma_{k}(A)=\sum_{|w|<k}(1-\frac{|w|}{k})a_{w}L_{w}.\]
Then the sequence $\Gamma_{k}(A)$ is strongly convergent to $A$.

By symmetry, for each $v$ in $\mathbb{F}_{n}^{+}$ we can define
an isometry $R_{v}$ by\[
R_{v}\xi_{w}=\xi_{wv},\quad w\in\mathbb{F}_{n}^{+},\]
and the map $v\to R_{v}$ gives an anti-representation (i.e. a multiplication
reversing representation) of $\mathbb{F}_{n}^{+}$, called the \emph{right
regular representation}. The isometries $R_{1},...,R_{n}$ also have
orthogonal ranges, and the free semigroup algebra they generate, denoted
by $\mathcal{R}_{n}$, is unitarily equivalent to $\mathcal{L}_{n}$.
It was shown in \cite{MR1665248} that $\mathcal{R}_{n}$ is the commutant
of $\mathcal{L}_{n}$.

A free semigroup algebra $\mathcal{S}$ is said to be of \emph{type
L} if it is algebraically isomorphic to $\mathcal{L}_{n}$. It was
shown in \cite{MR1823866} that if $\mathcal{S}$ is of type L, then
there is a completely isometric isomorphism $\Phi$ from $\mathcal{L}_{n}$
to $\mathcal{S}$ which takes the generators of $\mathcal{L}_{n}$
to the generators of $\mathcal{S}$. Moreover, $\Phi$ is a weak{*}-to-weak{*}
homeomorphism, and the inverse map $\Phi^{-1}$ is the dual of an
isometric isomorphism $\phi$ from the predual of $\mathcal{L}_{n}$
to the predual of $\mathcal{S}$.

For a free semigroup algebra $\mathcal{S}$, let $\mathcal{S}_{0}$
denote the weak-operator-closed ideal generated by $S_{1},...,S_{n}$.
Then either $\mathcal{S}_{0}=\mathcal{S}$, or $\mathcal{S}/\mathcal{S}_{0}\cong\mathbb{C}$.
In the latter case, the general structure theorem for free semigroup
algebras \cite{MR1823866} implies that $\mathcal{S}$ has a type
L part. If $\mathcal{S}_{0}=\mathcal{S}$, then $\mathcal{S}$ is
a von Neumann algebra.

The set of weak{*}-continuous linear functionals on $\mathcal{B}(\mathcal{H})$,
i.e. the predual, can be identified with the set of trace class operators
$\mathcal{C}^{1}(\mathcal{H})$, where $K$ in $\mathcal{C}^{1}(\mathcal{H})$
corresponds to the linear functional \[
T\to\operatorname{tr}(TK),\quad T\in\mathcal{B}(\mathcal{H}).\]
With this identification, the set of weak-operator-continuous linear
functionals on $\mathcal{B}(\mathcal{H})$ corresponds to the set
of finite rank operators. The predual of a weak{*}-closed subspace
$\mathcal{S}$ of $\mathcal{B}(\mathcal{H})$ can be identified with
the quotient space $\mathcal{C}^{1}(\mathcal{H})\diagup^{\perp}\mathcal{S}$,
where $^{\perp}\mathcal{S}$ denotes the set of elements in $\mathcal{C}^{1}(\mathcal{H})$
which annihilate $\mathcal{S}$, i.e. the preannihilator.

It was shown in \cite{MR1665248} that the weak{*} topology and the
weak operator topology coincide on $\mathcal{L}_{n}$. This means
that the coset of every weak{*}-continuous linear functional on $\mathcal{L}_{n}$
contains an element of finite rank. 

Let $\mathcal{S}$ be a free semigroup algebra on a Hilbert space
$\mathcal{H}$. A unit vector $x$ in $\mathcal{H}$ is said to be
\emph{wandering} for $\mathcal{S}$ if the set $\{S_{w}x:w\in\mathbb{F}_{n}^{+}\}$
is orthonormal. The following theorem from \cite{MR1823866} is integral
to our results.

\begin{thm}
\label{thm:amp-wand-vec}Let $\mathcal{S}$ be a type L free semigroup
algebra. Then for some $m\geq1$, the ampliation $\mathcal{S}^{(m)}$
has a wandering vector.
\end{thm}

Suppose that $\mathcal{S}$ is a type L free semigroup algebra, and
let $\pi_{0}$ be the weak-operator-continuous linear functional on
$\mathcal{S}$ such that $\pi_{0}$ annihilates $\mathcal{S}_{0}$
and $\pi_{0}(I)=1$. Then the coset in $\mathcal{C}^{1}(\mathcal{H})$
corresponding to $\pi_{0}$ contains an operator of finite rank, say
$m\geq1$. This $m$ corresponds to the $m$ in the statement of Theorem
\ref{thm:amp-wand-vec}. Since the restriction of $\mathcal{S}^{(m)}$
to the cyclic subspace generated by a wandering vector is unitarily
equivalent to $\mathcal{L}_{n}$, it follows that the weak{*} topology
and the weak operator topology agree on $\mathcal{S}$.

A subspace $\mathcal{S}$ of $\mathcal{B}(\mathcal{H})$ is said to
be \emph{reflexive} if $\mathcal{S}$ contains every operator $T$
in $\mathcal{B}(\mathcal{H})$ with the property that $Tx$ belongs
to $\mathcal{S}[x]$ for every $x$ in $\mathcal{H}$. This definition
of reflexivity was introduced by Loginov and Shulman \cite{MR0405124}.

The notion of hyper-reflexivity, which was introduced by Arveson \cite{MR0383098},
is a quantitative analogue of reflexivity. Let $d_{\mathcal{S}}$
denote the distance seminorm \[
d_{\mathcal{S}}(T)=\inf\{\|T-A\|:A\in\mathcal{S}\},\quad T\in\mathcal{B}(\mathcal{H}),\]
and define another seminorm $r_{\mathcal{S}}$ by \[
r_{\mathcal{S}}(T)=\sup\{|(Tx,y)|:\|x\|,\|y\|\leq1\mbox{ and }(Ax,y)=0\mbox{ for all }A\in\mathcal{S}\},\quad T\in\mathcal{B}(\mathcal{H}).\]
Then the reflexivity of $\mathcal{S}$ is equivalent to the condition
that $d_{\mathcal{S}}(T)=0$ if and only if $r_{\mathcal{S}}(T)=0$.

The equality $r_{\mathcal{S}}(T)\leq d_{\mathcal{S}}(T)$ always holds.
We say that $\mathcal{S}$ is \emph{hyper-reflexive} if there is a
constant $C>0$ such that $d_{\mathcal{S}}(T)\leq Cr_{\mathcal{S}}(T)$
for all $T$ in $\mathcal{B}(\mathcal{H})$. The smallest such $C$
is called the \emph{hyper-reflexivity constant} of $\mathcal{S}$.
Of course, hyper-reflexivity implies reflexivity.

Davidson \cite{MR882114} showed that the analytic Toeplitz algebra
is hyper-reflexive with hyper-reflexivity constant at most $19$.
Davidson and Pitts \cite{MR1665248} showed that for $n\geq2$, $\mathcal{L}_{n}$
is hyper-reflexive with hyper-reflexivity constant at most $51$.
This was later improved by Bercovici \cite{MR1641578}, who showed
that this hyper-reflexivity constant is at most $3$.

\section{\label{sec:toeplitz}The Non-commutative Toeplitz Operators}

The Toeplitz operators are precisely the operators $T$ in $\mathcal{B}(\ell^{2}(\mathbb{N}))$
which satisfy $S^{*}TS=T$, where $S$ is the unilateral shift. This
motivates the following definition, which was introduced by Popescu
\cite{MR1017331}.

\begin{defn}
Let $\mbox{S=}[S_{1}\enspace...\enspace S_{n}]$ be a row isometry.
We say that $T$ is an \emph{$S$-Toeplitz operator} if\[
S_{i}^{*}TS_{j}=\begin{cases}
T & \quad\mbox{if }i=j,\\
0 & \quad\mbox{otherwise},\end{cases}\]
and we let $\mathcal{T}_{S}$ denote the set of all $S$-Toeplitz
operators.
\end{defn}

If an $S$-Toeplitz operator $T$ is strictly positive, then by Theorem
4.3 of \cite{MR1017331}, it can be factored as $T=A^{*}A$, for some
$A$ in the commutant of the free semigroup algebra generated by $\mathcal{S}$.

Define row isometries $L$ and $R$ by $L=[L_{1}\enspace...\enspace L_{n}]$
and $R=[R_{1}\enspace...\enspace R_{n}]$. The size, $n$, will always
be clear from the context. In this section we will establish some
properties of the set $\mathcal{T}_{R}$ of $R$-Toeplitz operators
which we will need later. Note that since $\mathcal{L}_{n}$ is unitarily
equivalent to $\mathcal{R}_{n}$, the set $\mathcal{T}_{R}$ of $R$-Toeplitz
operators is unitarily equivalent to the set $\mathcal{T}_{L}$ of
$L$-Toeplitz operators. This means that any properties of $\mathcal{T}_{R}$
will correspond in an obvious way to properties of $\mathcal{T}_{L}$. 

The following Lemma is implied by Corollary 1.3 of \cite{MR2483229}.
Here we give a short direct proof.

\begin{lem}
\label{lem:set-of-r-toeplitz}The set $\mathcal{T}_{R}$ of $R$-Toeplitz
operators is precisely the weak{*} closure of the operator system
$\mathcal{L}_{n}^{*}+\mathcal{L}_{n}$.
\end{lem}
\begin{proof}
It is clear that the weak{*} closure of $\mathcal{L}_{n}^{*}+\mathcal{L}_{n}$
is contained in $\mathcal{T}_{R}$, since\[
R_{i}^{*}L_{w}R_{j}=R_{i}^{*}R_{j}L_{w}=\begin{cases}
L_{w} & \quad\mbox{if }i=j,\\
0 & \quad\mbox{otherwise}.\end{cases}\]
Suppose then that $T$ belongs to $\mathcal{T}_{R}$. It's clear that
$T^{*}$ also belongs to $\mathcal{T}_{R}$, and hence that the real
and imaginary parts of $T$ belong to $\mathcal{T}_{R}$. Since the
scalar operators also belong to $\mathcal{T}_{R}$, it follows that
we can write $T$ as a finite linear combination of strictly positive
operators in $\mathcal{T}_{R}$. Hence we may suppose that $T$ is
strictly positive.

By Theorem 4.3 of \cite{MR1017331}, we can write $T=A^{*}A$ for
some $A$ in $\mathcal{L}_{n}$. Note that $A^{*}\Gamma_{k}(A)$ belongs
to $\mathcal{L}_{n}+\mathcal{L}_{n}^{*}$ for $k\geq1$, where $\Gamma_{k}(A)$
denotes the $k$-th Cesaro sum of the Fourier series for $A$. The
sequence $\Gamma_{k}(A)$ is weak{*}-convergent to $A$, so it follows
that $A^{*}\Gamma_{k}(A)$ is weak{*}-convergent to $A^{*}A=T$, and
hence that $T$ belongs to the weak{*} closure of $\mathcal{L}_{n}+\mathcal{L}_{n}^{*}$.
\end{proof}

Note that based on the definition of the set $\mathcal{T}_{R}$ of
$R$-Toeplitz operators, Lemma \ref{lem:set-of-r-toeplitz} implies
that the weak{*} closure of $\mathcal{L}_{n}^{*}+\mathcal{L}_{n}$
is closed in the weak operator topology. 

\begin{lem}
\label{lem:factor-r-toeplitz}For $n\geq2$, every $R$-Toeplitz operator
$T$ can be factored as $T=B^{*}C$ for some $B$ and $C$ in $\mathcal{L}_{n}$.
Moreover, $B$ and $C$ can be taken to be bounded below.
\end{lem}
\begin{proof}
As in the proof of Lemma \ref{lem:set-of-r-toeplitz}, we can write
$T$ as a finite linear combination of strictly positive $R$-Toeplitz
operators, say $T=\sum_{i=1}^{m}c_{i}T_{i}$ for some $c_{1},...,c_{m}$
in $\mathbb{C}$ and strictly positive $T_{1},...,T_{m}$ in $\mathcal{T}_{R}$.
By Theorem 4.3 of \cite{MR1017331}, we can factor each $T_{i}$ as
$T_{i}=A_{i}^{*}A_{i}$ for some $A_{i}$ in $\mathcal{L}_{n}$. Set
$B=\sum_{i=1}^{m}L_{1^{i}2}A_{i}$ and $C=\sum_{i=1}^{m}c_{i}L_{1^{i}2}A_{i}$.
Then $B$ and $C$ both belong to $\mathcal{L}_{n}$ and $T=B^{*}C$.

To see that $B$ and $C$ can be taken to be bounded below, take $B'=B+L_{1^{m+1}2}$
and $C'=C+L_{1^{m+2}2}$. Then $B'$ and $C'$ both belong to $\mathcal{L}_{n}$.
Since the isometries $L_{12},...,L_{1^{m+2}2}$ have pairwise orthogonal
ranges, $B'$ and $C'$ are bounded below, and $T=(B')^{*}C'$. 
\end{proof}

Lemma \ref{lem:factor-r-toeplitz} provides another characterization
of the $R$-Toeplitz operators for $n\geq2$.

\begin{cor}
\label{cor:char-r-toeplitz}For $n\geq2$, the set $\mathcal{T}_{R}$
of $R$-Toeplitz operators is precisely $\mathcal{L}_{n}^{*}\mathcal{L}_{n}=\{B^{*}C:B,C\in\mathcal{L}_{n}\}$.
\end{cor}

Popescu \cite{MR2483229} showed that every $R$-Toeplitz operator
$T$ has a Fourier series\begin{eqnarray*}
T & \sim & \sum_{w\in\mathbb{F}_{n}^{+}}a_{w}L_{w}+\sum_{w\in\mathbb{F}_{n}^{+}\setminus\{\varnothing\}}\overline{b_{w}}L_{w}^{*},\end{eqnarray*}
which is a formal power series with coefficients in $\mathcal{L}_{n}$
and $\mathcal{L}_{n}^{*}$. This completely determines $T$ in the
sense that for every word $u$ in $\mathbb{F}_{n}^{+}$, \[
T\xi_{u}=\sum_{w\in\mathbb{F}_{n}^{+}}a_{w}L_{w}\xi_{u}+\sum_{w\in\mathbb{F}_{n}^{+}\setminus\{\varnothing\}}\overline{b_{w}}L_{w}^{*}\xi_{u}.\]

Let $\mathcal{S}$ be a type L free semigroup algebra. We know that
the canonical map from $\mathcal{L}_{n}$ to $\mathcal{S}$ is a complete
isometry and a weak{*}-to-weak{*} homeomorphism. Our goal for the
remainder of this section is to show that this map extends in a natural
way to a map from the weak{*} closure of $\mathcal{L}_{n}+\mathcal{L}_{n}^{*}$
(i.e. from the set $\mathcal{T}_{R}$ of $R$-Toeplitz operators)
to the weak{*} closure of $\mathcal{S}+\mathcal{S}^{*}$, and that
this extension is also a complete isometry and a weak{*}-to-weak{*}
homeomorphism.

\begin{lem}
\label{lem:canon-isom-to-isom}Let $\mathcal{S}$ be a type L free
semigroup algebra with $n\geq2$ generators , and let $\Phi$ be the
canonical map from $\mathcal{L}_{n}$ to $\mathcal{S}$. Then $\Phi^{-1}$
maps isometries in $\mathcal{S}$ to isometries in $\mathcal{L}_{n}$.
\end{lem}
\begin{proof}
By Theorem \ref{thm:amp-wand-vec}, $\mathcal{S}^{(m)}$ has a wandering
vector $w$ for some $m$, and the restriction of $\mathcal{S}^{(m)}$
to $\mathcal{S}^{(m)}[w]$ is unitarily equivalent to $\mathcal{L}_{n}$.
The map $\Phi^{-1}$ from $\mathcal{S}$ to $\mathcal{L}_{n}$ is
given by taking $\mathcal{S}$ to $\mathcal{S}^{(m)}$ , restricting
to $\mathcal{S}^{(m)}[w]$, and applying this equivalence. If $G$
is an isometry in $\mathcal{S}$, then $G^{(m)}$ is an isometry in
$\mathcal{S}^{(m)}$, and so clearly the restriction of $G^{(m)}$
to $\mathcal{S}^{(m)}[w]$ is an isometry.
\end{proof}

\begin{thm}
\label{thm:exten-canon}Let $\mathcal{S}$ be a type L free semigroup
algebra with $n\geq2$ generators on a Hilbert space $\mathcal{H}$.
Then the canonical map $\Phi$ from $\mathcal{L}_{n}$ to $\mathcal{S}$
extends to a completely isometric weak{*}-to-weak{*} homeomorphism
from the weak{*} closure of $\mathcal{L}_{n}+\mathcal{L}_{n}^{*}$
to the weak{*} closure of $\mathcal{S}+\mathcal{S}^{*}$. 
\end{thm}
\begin{proof}
Applying Arveson's extension theorem \cite{MR0253059} gives a completely
positive map $\Psi$ from $C^{*}(\mathcal{L}_{n})$ to $\mathcal{B}(\mathcal{H})$
which extends $\Phi$. Since $\Psi$ extends $\Phi$, we have $\|\Psi\|=\|\Psi(I)\|=\|\Phi(I)\|=1$.
Let $\mathcal{Z}=\{A\in C^{*}(\mathcal{L}_{n}):\Psi(A)^{*}\Psi(A)=\Psi(A^{*}A)\}$.
By \cite{MR0355615}, we have\[
\mathcal{Z}=\{A\in C^{*}(\mathcal{L}_{n}):\Psi(B)\Psi(A)=\Psi(BA)\mbox{ for all }B\mbox{ in }C^{*}(\mathcal{L}_{n})\}.\]
By Theorem 4.1 of \cite{MR1823866}, $\Phi$ maps isometries in $\mathcal{L}_{n}$
to isometries in $\mathcal{S}$, so every isometry in $\mathcal{L}_{n}$
belongs to $\mathcal{Z}$. Since, by Theorem 4.5 of \cite{MR1823866},
every element in $\mathcal{L}_{n}$ can be written as a finite linear
combination of isometries in $\mathcal{L}_{n}$, this implies that
$\mathcal{Z}$ contains all of $\mathcal{L}_{n}$. Hence for $A$
in $\mathcal{L}_{n}$, $\Psi(TA)=\Psi(T)\Psi(A)$ for all $T$ in
$C^{*}(\mathcal{L}_{n})$. Note that by Corollary \ref{cor:char-r-toeplitz},
$C^{*}(\mathcal{L}_{n})$ contains $\mathcal{T}_{R}$. For the remainder
of the proof, we restrict $\Psi$ to $\mathcal{T}_{R}$.

Let $T$ be a self-adjoint element in $\mathcal{T}_{R}$ such that
$\Psi(T)=0$. For sufficiently large $\lambda>0$, $T+\lambda I$
is strictly positive, so by Theorem 4.3 of \cite{MR1017331}, we can
write $T+\lambda I=B^{*}B$ for some $B$ in $\mathcal{L}_{n}$. Let
$V=\lambda^{-1/2}B$. Then\begin{eqnarray*}
\Phi(V)^{*}\Phi(V)-I & = & \Psi(V^{*}V-I)\\
 & = & \Psi(\lambda^{-1}B^{*}B-I)\\
 & = & \Psi(\lambda^{-1}(T+\lambda I)-I)\\
 & = & \lambda^{-1}\Psi(T)\\
 & = & 0,\end{eqnarray*}
which shows that $\Phi(V)$ is an isometry in $\mathcal{S}$. By Lemma
\ref{lem:canon-isom-to-isom}, this implies that $V$ is an isometry
in $\mathcal{L}_{n}$. Hence \[
T=\lambda(V^{*}V-I)=0.\]
Since, for arbitrary $T$ in $\mathcal{T}_{R}$, $\operatorname{re}(T)$
and $\operatorname{im}(T)$ are self-adjoint, and since \[
\Psi(T)=\Psi(\operatorname{re}(T)+\operatorname{im}(T))=\operatorname{re}(\Psi(T))+\operatorname{im}(\Psi(T))=0\]
 if and only if $\psi(\operatorname{re}(T))=0$ and $\psi(\operatorname{im}(T))=0$,
it follows that $\Psi$ is injective.

Arguing exactly as above, the canonical map $\Phi^{-1}$ from $\mathcal{S}$
to $\mathcal{L}_{n}$ also has a completely positive extension $\Omega$
from $C^{*}(\mathcal{S})$ to $\mathcal{B}(\mathcal{F}_{n})$, and
for $G$ in $\mathcal{S}$, $\Omega(HG)=\Omega(H)\Omega(G)$ for all
$H$ in $C^{*}(\mathcal{S})$. Since $\Omega$ extends $\Phi^{-1}$,
we have $\|\Omega\|=\|\Omega(I)\|=\|\Phi^{-1}(I)\|=1$. For the remainder
of the proof we restrict $\Omega$ to the intersection of $C^{*}(\mathcal{S})$
and the range of $\Psi$.

Note that the range of $\Psi$ is contained in the weak{*} closure
of $\mathcal{S}+\mathcal{S}^{*}$. Indeed, by Lemma \ref{cor:char-r-toeplitz},
every element in the range of $\Psi$ can be written as $\Psi(B^{*}C)=\Psi(B^{*})\Psi(C)=\Phi(B)^{*}\Phi(C)$
for some $B$ and $C$ in $\mathcal{L}_{n}$. The sequence $\Gamma_{k}(C)$
is weak operator convergent to $C$, so by the weak operator continuity
of $\Phi$, the sequence $\Phi(\Gamma_{k}(C))$ is weak operator convergent
to $\Phi(C)$, and hence the sequence $\Phi(B)^{*}\Phi(\Gamma_{k}(C))$
is weak operator convergent to $\Phi(B)^{*}\Phi(C)$, which implies
that $\Phi(B)^{*}\Phi(C)$ is contained in the weak{*} closure of
$\mathcal{S}+\mathcal{S}^{*}$.

We claim that $\Omega(\Psi(T))=T$ for all $T$ in $\mathcal{T}_{R}$.
Indeed, apply Lemma \ref{cor:char-r-toeplitz} to write $T=B^{*}C$
for some $B$ and $C$ in $\mathcal{L}_{n}$, and let $G=\Phi(B)$
and $H=\Phi(C)$. Then we have\begin{eqnarray*}
\Psi(T) & = & \Psi(B^{*}C)\\
 & = & \Phi(B)^{*}\Phi(C)\\
 & = & G^{*}H,\end{eqnarray*}
which gives\begin{eqnarray*}
\Omega(\Phi(T)) & = & \Omega(G^{*}H)\\
 & = & (\Phi^{-1}(G))^{*}\Phi^{-1}(H)\\
 & = & B^{*}C\\
 & = & T.\end{eqnarray*}
Then \[
\frac{\|T\|}{\|\Psi(T)\|}=\frac{\|\Omega(\Psi(T))\|}{\|\Psi(T)\|}\leq1,\]
which gives\[
\|T\|\leq\|\Psi(T)\|\leq\|T\|,\]
and shows that $\Psi$ maps $\mathcal{T}_{R}$ isometrically onto
its range.

We now show that $\Psi$ is weak{*}-to-weak{*} continuous. Since the
predual of $\mathcal{T}_{R}$ is separable, by an application of the
Krein-Smulian theorem it suffices to show that if $T_{n}$ is a sequence
in $\mathcal{T}_{R}$ which is weak{*}-convergent to zero, then $\Psi(T_{n})$
is weak{*} convergent to zero.

Let $\mathcal{A}=\{A\oplus\Phi(A):A\in\mathcal{L}_{n}\}$, and note
that $\mathcal{A}$ is the free semigroup algebra generated by the
isometries $L_{1}\oplus S_{1},...,L_{n}\oplus S_{n}$. Fix $u$ in
$\mathcal{H}$. By Theorem 1.6 of \cite{MR1823866}, there exists
a vector $x$ in $\mathcal{F}_{n}$ such that the restriction of $\mathcal{A}$
to $\mathcal{W}=\mathcal{A}[x\oplus u]$ is unitarily equivalent to
$\mathcal{L}_{n}$. Letting $P$ denote the projection of $\mathcal{F}\oplus\mathcal{H}$
onto $\mathcal{W}$, and letting $\mathcal{K}$ denote the weak{*}
closure of the restriction of $P(\mathcal{A}+\mathcal{A}^{*})P$ to
$\mathcal{W}$, it follows that $\mathcal{K}$ is unitarily equivalent
to $\mathcal{T}_{R}$. By Lemma \ref{cor:char-r-toeplitz}, every
element of $\mathcal{K}$ can be written as the restriction to $\mathcal{W}$
of an element of the form\[
P(B^{*}\oplus\Phi(B)^{*})(C\oplus\Phi(C))P=P(B^{*}C\oplus\Psi(B^{*}C))P.\]
Hence $\mathcal{K}$ is the restriction to $\mathcal{W}$ of $\{T\oplus\psi(T):T\in\mathcal{T}_{R}\}$.

If $T_{n}$ is weak{*} convergent to zero in $\mathcal{T}_{R}$, the
unitary equivalence between $\mathcal{T}_{R}$ and $\mathcal{K}$
implies the restriction of the sequence $T_{n}\oplus\Psi(T_{n})$
to $\mathcal{W}$ is weak{*}-convergent to zero in $\mathcal{K}$.
Hence\[
((T_{n}\oplus\Psi(T_{n}))(x\oplus u),x\oplus u)=(T_{n}x,x)+(\Psi(T_{n})u,u)\to0,\]
and since $(T_{n}x,x)\to0$, this implies that $(\Psi(T_{n})u,u)\to0$.
Since $u$ was chosen arbitrarily, we deduce that $(\Psi(T_{n})u,u)\to0$
for all $u$ in $\mathcal{H}$. By the polarization identity, we get
that $\Psi(T_{n})$ is weak operator convergent to zero. By the uniform
boundedness principle, the sequence $\Psi(T_{n})$ is bounded. It
follows that $\Psi(T_{n})$ is weak{*} convergent to zero. We therefore
conclude that $\Psi$ is weak{*} continuous.

It now follows by another application of the Krein-Smulian theorem
that $\Psi$ has weak{*} closed range, and that $\Psi$ is a weak{*}-to-weak{*}
homeomorphism onto its range. But it's clear that the range of $\Psi$
is weak{*} dense in the weak{*} closure of $\mathcal{S}+\mathcal{S}^{*}$,
so $\Psi$ maps $\mathcal{T}_{R}$ weak{*}-to-weak{*} homeomorphically
onto the weak{*} closure of $\mathcal{S}+\mathcal{S}^{*}$. From above,
$\Psi$ is a completely positive isometry, with completely positive
inverse $\Omega$. Hence $\Psi$ is completely isometric.
\end{proof}

\section{\label{sec:wandering}Wandering Vectors}

Let $\mathcal{S}$ be a weak{*}-closed subspace of $\mathcal{B}(\mathcal{H})$,
and let $x$ and $y$ be vectors in $\mathcal{H}$. Then $[x\otimes y]_{\mathcal{S}}$
denotes the weak-operator-continuous linear functional on $\mathcal{S}$
which is given by the coset of the rank one tensor $x\otimes y$.
In other words, \[
(A,[x\otimes y]_{\mathcal{S}})=(Ax,y),\quad A\in\mathcal{S}.\]

\begin{defn}
A weak{*}-closed subspace $\mathcal{S}$ of $\mathcal{B}(\mathcal{H})$
is said to have property $\mathbb{A}_{1}(1)$ if, for every weak{*}-continuous
linear functional $\pi$ on $\mathcal{S}$ and every $\varepsilon>0$,
there are vectors $x$ and $y$ in $\mathcal{H}$ with $\|x\|\|y\|<(1+\epsilon)\|\pi\|$
such that $\pi(A)=(Ax,y)$ for all $A$ in $\mathcal{S}$.
\end{defn}

It was shown in \cite{MR1665248} that $\mathcal{L}_{n}$ has property
$\mathbb{A}_{1}(1)$. A result of Bercovici \cite{MR960951} implies
that a singly generated type L free semigroup algebra has property
$\mathbb{A}_{1}(1)$. In this section, we will use dual algebra techniques
to show that every type L free semigroup algebra with $n\geq2$ generators
has property $\mathbb{A}_{1}(1)$. From this result, it will follow
easily that every type L free semigroup algebra has a wandering vector.

For the remainder of this section we fix a type L free semigroup algebra
$\mathcal{S}$ with $n\geq2$ generators acting on a Hilbert space
$\mathcal{H}$. The general outline of our approach is as follows.
Let $\pi$ be a weak{*}-continuous linear functional on $\mathcal{S}$.
We will show that we can construct convergent sequences $(x_{k})$
and $(y_{k}$) such that\[
\lim_{k\to\infty}\|\pi-[x_{k}\otimes y_{k}]_{\mathcal{S}}\|=0.\]
This will then give $\pi=[x\otimes y]_{\mathcal{S}}$, where $x=\lim_{k}x_{k}$
and $y=\lim_{k}y_{k}$.

The following idea will allow us to iteratively construct the sequences
$(x_{k})$ and $(y_{k})$. Fix $x_{k}$ and $y_{k}$. Suppose we can
find vectors $x'$ and $y'$ such that

\begin{enumerate}
\item $[x'\otimes y']_{\mathcal{S}}$ approximates the error $\pi-[x_{k}\otimes y_{k}]_{\mathcal{S}}$
arbitrarily closely,
\item $\|[x_{k}\otimes y']_{\mathcal{S}}\|$ and $\|[x'\otimes y_{k}]_{\mathcal{S}}\|$
are arbitrarily small,
\item $\|x'\|$ and $\|y'\|$ are arbitrarily close to $\|\pi-[x_{k}\otimes y_{k}]_{\mathcal{S}}\|.$
\end{enumerate}
Set $x_{k+1}=x_{k}+x'$ and $y_{k+1}=y_{k}+y'$. Then\[
\|\pi-[x_{k+1}\otimes y_{k+1}]_{\mathcal{S}}\|\leq\|\pi-[x_{k}\otimes y_{k}]_{\mathcal{S}}-[x'\otimes y']_{\mathcal{S}}\|+\|[x'\otimes y_{k}]_{\mathcal{S}}\|+\|[x_{k}\otimes y']_{\mathcal{S}}\|,\]
so $[x_{k+1}\otimes y_{k+1}]_{\mathcal{S}}$ is an arbitrarily good
approximation to $\pi$, and the sequences $(x_{k})$ and $(y_{k})$
can be made Cauchy. Of course, the main difficulty will be in showing
that it is possible to find $x'$ and $y'$ as above.

\begin{defn}
An operator $X:\mathcal{F}_{n}\to\mathcal{H}$ is said to \emph{intertwine}
$\mathcal{L}_{n}$ and $\mathcal{S}$ if $XL_{i}=S_{i}X$ for $1\leq i\leq n$.
\end{defn}

Let $\overline{x}=(x_{1},...,x_{m})$ be a wandering vector for $\mathcal{S}^{(m)}$.
We know that the restriction of $\mathcal{S}^{(m)}$ to $\mathcal{S}^{(m)}[\overline{x}]$
is unitarily equivalent to $\mathcal{L}_{n}$. Let $X:\mathcal{F}_{n}\to\mathcal{H}$
denote the map which follows this equivalence with the projection
onto the first coordinate. Then $X$ intertwines $\mathcal{L}_{n}$
and $\mathcal{S}$. It was shown in \cite{MR2139108} that every vector
in $\mathcal{H}$ is in the range of some intertwining operator of
this form.

The following result shows that every intertwining operator gives
rise to an $L$-Toeplitz operator. This allows us to use the results
of section \ref{sec:toeplitz} to work with intertwining operators.

\begin{lem}
\label{lem:intertwine-gives-toeplitz}Suppose $X:\mathcal{F}_{n}\to\mathcal{H}$
intertwines $\mathcal{L}_{n}$ and $\mathcal{S}$. Then $X^{*}X$
is an $L$-Toeplitz operator.
\end{lem}
\begin{proof}
This follows immediately from the identity\begin{eqnarray*}
L_{i}^{*}X^{*}XL_{j} & = & X^{*}S_{i}^{*}S_{j}X\\
 & = & \begin{cases}
X^{*}X & \quad\mbox{if }i=j,\\
0 & \quad\mbox{otherwise}.\end{cases}\end{eqnarray*}

\end{proof}

We require several technical results about $L$-Toeplitz operators.

\begin{lem}
\label{lem:short-ineq}Let $T$ be an $L$-Toeplitz operator with
Fourier series\[
T\sim\sum_{w\in\mathbb{F}_{n}^{+}}a_{w}R_{w}+\sum_{w\in\mathbb{F}_{n}^{+}\setminus\{\varnothing\}}\overline{b_{w}}R_{w}^{*}.\]
Then for any word $u$ in $\mathbb{F}_{n}^{+}$,\[
\left\Vert \sum_{w\in\mathbb{F}_{n}^{+}\setminus\{\varnothing\}}\overline{b_{w}}R_{w}^{*}\xi_{u}\right\Vert \leq\left\Vert \sum_{w\in\mathbb{F}_{n}^{+}\setminus\{\varnothing\}}b_{w}R_{w}\xi_{u}\right\Vert .\]

\end{lem}
\begin{proof}
We have\begin{eqnarray*}
\left\Vert \sum_{w\in\mathbb{F}_{n}^{+}\backslash\{\varnothing\}}\overline{b_{w}}R_{w}^{*}\xi_{u}\right\Vert ^{2} & = & \sum_{\substack{w\in\mathbb{F}_{n}^{+}\backslash\{\varnothing\}\\
w=w'u}
}\left|b_{w}\right|^{2}\\
 & \leq & \sum_{w\in\mathbb{F}_{n}^{+}\backslash\{\varnothing\}}\left|b_{w}\right|^{2}\\
 & = & \left\Vert \sum_{w\in\mathbb{F}_{n}^{+}\backslash\{\varnothing\}}b_{w}R_{w}\xi_{u}\right\Vert ^{2}.\end{eqnarray*}

\end{proof}

\begin{lem}
\label{lem:flattening-word}Let $T$ be an $L$-Toeplitz operator
with Fourier series\[
T\sim\sum_{w\in\mathbb{F}_{n}^{+}}a_{w}R_{w}+\sum_{w\in\mathbb{F}_{n}^{+}\setminus\{\varnothing\}}\overline{b_{w}}R_{w}^{*}.\]
Then given $p\geq1$ and $\epsilon>0$, there is a word $v$ in $\mathbb{F}_{n}^{+}$
such that \[
\|R_{v}^{*}TR_{v}\xi_{u}-a_{\varnothing}\xi_{u}\|<\epsilon\]
 for any word $u\in\mathbb{F}_{n}^{p}$.
\end{lem}
\begin{proof}
For $k\geq1$, let $v_{k}$ be the word $v_{k}=12^{k}$. Then for
any word $w$ in $\mathbb{F}_{n}^{+}$,\begin{eqnarray*}
R_{v_{k}}^{*}R_{w}R_{v_{k}} & = & \begin{cases}
I & \mbox{if }w=\varnothing,\\
R_{v_{k}w'} & \mbox{if }w=w'v_{k}\mbox{ for }w'\in\mathbb{F}_{n}^{+},\\
0 & \mbox{otherwise}.\end{cases}\end{eqnarray*}
This implies the Fourier series for $R_{v_{k}}^{*}TR_{v_{k}}$ is
given by\[
R_{v_{k}}^{*}TR_{v_{k}}\sim a_{\varnothing}I+\sum_{\substack{w\in\mathbb{F}_{n}^{+}\backslash\{\varnothing\}\\
w=w'v_{k}}
}a_{w}R_{v_{k}w'}+\sum_{\substack{w\in\mathbb{F}_{n}^{+}\backslash\{\varnothing\}\\
w=w'v_{k}}
}\overline{b_{w}}R_{v_{k}w'}^{*}.\]
Hence for $u$ in $\mathbb{F}_{n}^{+}$,\[
R_{v_{k}}^{*}TR_{v_{k}}\xi_{u}\sim a_{\varnothing}\xi_{u}+\sum_{\substack{w\in\mathbb{F}_{n}^{+}\backslash\{\varnothing\}\\
w=w'v_{k}}
}a_{w}R_{v_{k}w'}\xi_{u}+\sum_{\substack{w\in\mathbb{F}_{n}^{+}\backslash\{\varnothing\}\\
w=w'v_{k}}
}\overline{b_{w}}R_{v_{k}w'}^{*}\xi_{u}.\]
This gives\begin{eqnarray*}
\|R_{v_{k}}^{*}TR_{v_{k}}\xi_{u}-a_{\varnothing}\xi_{u}\| & = & \left\Vert \sum_{\substack{w\in\mathbb{F}_{n}^{+}\backslash\{\varnothing\}\\
w=w'v_{k}}
}a_{w}R_{v_{k}w'}\xi_{u}+\sum_{\substack{w\in\mathbb{F}_{n}^{+}\backslash\{\varnothing\}\\
w=w'v_{k}}
}\overline{b_{w}}R_{v_{k}w'}^{*}\xi_{u}\right\Vert \\
 & \leq & \left\Vert \sum_{\substack{w\in\mathbb{F}_{n}^{+}\backslash\{\varnothing\}\\
w=w'v_{k}}
}a_{w}R_{v_{k}w'}\xi_{u}\right\Vert +\left\Vert \sum_{\substack{w\in\mathbb{F}_{n}^{+}\backslash\{\varnothing\}\\
w=w'v_{k}}
}\overline{b_{w}}R_{v_{k}w'}^{*}\xi_{u}\right\Vert \\
 & \leq & \left\Vert \sum_{\substack{w\in\mathbb{F}_{n}^{+}\backslash\{\varnothing\}\\
w=w'v_{k}}
}a_{w}R_{v_{k}w'}\xi_{u}\right\Vert +\left\Vert \sum_{\substack{w\in\mathbb{F}_{n}^{+}\backslash\{\varnothing\}\\
w=w'v_{k}}
}b_{w}R_{v_{k}w'}\xi_{u}\right\Vert ,\end{eqnarray*}
where the last inequality follows from Lemma \ref{lem:short-ineq}.
Now\[
\left\Vert \sum_{\substack{w\in\mathbb{F}_{n}^{+}\backslash\{\varnothing\}\\
w=w'v_{k}}
}a_{w}R_{v_{k}w'}\xi_{u}\right\Vert ^{2}=\sum_{\substack{w\in\mathbb{F}_{n}^{+}\backslash\{\varnothing\}\\
w=w'v_{k}}
}\left|a_{w}\right|^{2}=\|R_{v_{k}}^{*}T\xi_{\varnothing}\|^{2},\]
and similarly,\[
\left\Vert \sum_{\substack{w\in\mathbb{F}_{n}^{+}\backslash\{\varnothing\}\\
w=w'v_{k}}
}b_{w}R_{v_{k}w'}\xi_{u}\right\Vert ^{2}=\sum_{\substack{w\in\mathbb{F}_{n}^{+}\backslash\{\varnothing\}\\
w=w'v_{k}}
}\left|b_{w}\right|^{2}=\|R_{v_{k}}^{*}T^{*}\xi_{\varnothing}\|^{2},\]
Hence the result follows from the fact that for all $\xi$ in $\mathcal{F}_{n}$,
$\|R_{v}^{*}\xi\|\to0$ as $\left|v\right|\to\infty$.
\end{proof}

Recall that $\phi:(\mathcal{L}_{n})_{*}\to\mathcal{S}_{*}$ is the
predual of the map $\Phi^{-1}:\mathcal{S}\to\mathcal{L}_{n}$.

\begin{lem}
\label{lem:word-close-to-wand}Let $X:\mathcal{F}_{n}\to\mathcal{H}$
be an intertwining operator, and let $x=X\xi_{\varnothing}$. Then
given $p\geq1$ and $\epsilon>0$, there exists a word $v$ in $\mathbb{F}_{n}^{+}$
such that \[
\left\Vert [S_{u_{1}v}x\otimes S_{u_{2}v}x]_{\mathcal{S}}-\|x\|^{2}\phi([\xi_{u_{1}}\otimes\xi_{u_{2}}]_{\mathcal{L}_{n}})\right\Vert <\epsilon\]
 for all words $u_{1}$ and $u_{2}$ in $\mathbb{F}_{n}^{p}$.
\end{lem}
\begin{proof}
By scaling $X$ if necessary, we can suppose that $\|x\|=1$. Let
$T=X^{*}X$. Then $T$ is an $L$-Toeplitz operator by Lemma \ref{lem:intertwine-gives-toeplitz}.
Writing the Fourier series for $T$ as\[
T\sim\sum_{w\in\mathbb{F}_{n}^{+}}a_{w}R_{w}+\sum_{w\in\mathbb{F}_{n}^{+}\setminus\{0\}}\overline{b_{w}}R_{w}^{*},\]
it follows that $a_{\varnothing}=\|x\|^{2}=1$. Hence by Lemma \ref{lem:flattening-word},
there exists a word $v$ in $\mathbb{F}_{n}^{+}$ such that $\|R_{v}^{*}TR_{v}\xi_{u_{2}}-\xi_{u_{2}}\|<\epsilon$
for any word $u_{2}$ in $\mathbb{F}_{n}^{p}$. Then for $A$ in $\mathcal{S}$,
\begin{eqnarray*}
(A,[S_{u_{1}v}x\otimes S_{u_{2}v}x]_{\mathcal{S}}) & = & (AS_{u_{1}v}x,S_{u_{2}v}x)\\
 & = & (AS_{u_{1}v}X\xi_{\varnothing},S_{u_{2}v}X\xi_{\varnothing})\\
 & = & (X\Phi^{-1}(A)L_{u_{1}v}\xi_{\varnothing},XL_{u_{2}v}\xi_{\varnothing})\\
 & = & (X\Phi^{-1}(A)R_{v}\xi_{u_{1}},XR_{v}\xi_{u_{2}})\\
 & = & (XR_{v}\Phi^{-1}(A)\xi_{u_{1}},XR_{v}\xi_{u_{2}})\\
 & = & (\Phi^{-1}(A)\xi_{u_{1}},R_{v}^{*}TR_{v}\xi_{u_{2}})\end{eqnarray*}
for all words $u_{1}$ and $u_{2}$ in $\mathbb{F}_{n}^{p}$. This
gives\begin{eqnarray*}
\left|(A,[S_{u_{1}v}x\otimes S_{u_{2}v}x]_{\mathcal{S}}-[\xi_{u_{1}}\otimes\xi_{u_{2}}]_{\mathcal{L}_{n}})\right| & = & \left|(\Phi^{-1}(A)\xi_{u_{1}},R_{v}^{*}TR_{v}\xi_{u_{2}}-\xi_{u_{2}})\right|\\
 & \leq & \|\Phi^{-1}(A)\xi_{u_{1}}\|\|R_{v}^{*}TR_{v}\xi_{u_{2}}-\xi_{u_{2}}\|\\
 & < & \epsilon\|A\|.\end{eqnarray*}
Therefore,\[
\|[S_{u_{1}v}x\otimes S_{u_{2}v}x]_{\mathcal{S}}-\phi([\xi_{u_{1}}\otimes\xi_{u_{2}}]_{\mathcal{L}_{n}})\|<\epsilon.\qedhere\]

\end{proof}

\begin{lem}
\label{lem:key-lemma-2}Let $M\geq2$ be minimal such that $\mathcal{S}^{(M)}$
has a wandering vector $\overline{w}=(w_{1},...,w_{M})$. Then given
$\epsilon\in(0,1)$ there exists a unit vector $\overline{x}=(x_{1},...,x_{M})$
in $\mathcal{S}^{(M)}[\overline{w}]$ such that $x_{1}=X\xi_{\varnothing}$
for some intertwining operator $X:\mathcal{F}_{n}\to\mathcal{H}$,
and $\|x_{1}\|>1-\epsilon$.
\end{lem}
\begin{proof}
Let $P$ denote the projection map from $\mathcal{S}^{(M)}[\overline{w}]$
to $\mathcal{H}^{(M-1)}$ which takes $\overline{x}=(x_{1},...,x_{M})$
to $(x_{2},...,x_{M})$. Then $P$ intertwines the restriction of
$\mathcal{S}^{(M)}$ to $\mathcal{S}^{(M)}[\overline{w}]$ and $\mathcal{S}^{(M-1)}$.
The restriction of $\mathcal{S}^{(M)}$ to $\mathcal{S}^{(M)}[\overline{w}]$
is unitarily equivalent to $\mathcal{L}_{n}$. Let $U$ be a unitary
implementing this equivalence. Then setting $Y=PU$, $Y$ intertwines
$\mathcal{L}_{n}$ and $S^{(M-1)}$. Suppose that for all $\overline{x}$
in $\mathcal{S}^{(M)}[\overline{w}]$, $\|x_{1}\|\leq(1-\epsilon)\|\overline{x}\|$.
Then\[
\|\overline{x}\|^{2}=\sum_{i=1}^{M}\|x_{i}\|^{2}\leq(1-\epsilon)^{2}\|\overline{x}\|^{2}+\sum_{i=2}^{M}\|x_{i}\|^{2},\]
which gives\[
\sum_{i=2}^{mM}\|x_{i}\|^{2}\geq(1-(1-\epsilon)^{2})\|\overline{x}\|^{2},\]
implying that $P$ is bounded below, and hence that $Y$ is bounded
below. By \cite[Theorem 2.8]{MR2139108}, this implies that the range
of $Y$ is a wandering subspace for $\mathcal{S}^{(M-1)}$, contradicting
the minimality of $M$. Hence there must be some unit vector $\overline{x}$
in $\mathcal{S}^{(M)}[\overline{w}]$ such that $\|x_{1}\|>1-\epsilon$.

Let $Q$ denote the projection map from $\mathcal{S}^{(M)}[\overline{w}]$
to $\mathcal{H}$ which takes $\overline{y}=(y_{1},...,y_{M})$ to
$\overline{y}_{1}$, and let $Z=QU$. Note that $x_{1}$ is contained
in the range of $Z$. For every $R$ in $\mathcal{R}_{n}$, the operator
$ZR$ intertwines $\mathcal{L}_{n}$ and $\mathcal{S}$. Moreover,
since the set of vectors $\{R\xi_{\varnothing}:R\in\mathcal{R}_{n}\}$
is dense in $\mathcal{F}_{n}$, the set $\{ZR\xi_{\varnothing}:R\in\mathcal{R}_{n}\}$
is dense in the closure of the range of $Z$. It follows that we can
choose the vector $\overline{x}$ as above such that $x_{1}=X\xi_{\varnothing}$
for some intertwining operator $X:\mathcal{F}_{n}\to\mathcal{H}$.
\end{proof}

Let $M\geq1$ be minimal such that the ampliation $\mathcal{S}^{(M)}$
has a wandering vector $\overline{x}$. Such $M$ exists by Theorem
\ref{thm:amp-wand-vec}. Then $\mathcal{H}^{(M)}$ contains an infinite
family of pairwise orthogonal subspaces $\mathcal{W}_{k}$, for $k\geq1$,
which are wandering for $\mathcal{S}^{(M)}$. For example, we can
take $\mathcal{W}_{k}=\mathcal{S}^{(m)}[S_{v_{k}}^{(m)}\overline{x}]$,
where $v_{k}=12^{k}$. For $k\geq1$, let $\mathcal{M}_{k}$ denote
the linear manifold in $\mathcal{H}$ given by\[
\mathcal{M}_{k}=\{z\in\mathcal{H}:z=\overline{z}_{1}\mbox{ for some }\overline{z}=(\overline{z}_{1},...,\overline{z}_{M})\mbox{ in }\mathcal{W}_{k}\}.\]
Let $\mathcal{W}$ denote the algebraic span of the $\mathcal{W}_{k}$,
and let $\mathcal{M}$ denote the algebraic span of the $\mathcal{M}_{k}$.

\begin{lem}
\label{lem:quasi-orthog}Given $h_{1},...,h_{q}$ in $\mathcal{M}$
and $\epsilon>0$, there exists a unit vector $y$ in $\mathcal{M}$
such that $y=Y\xi_{\varnothing}$ for some intertwining operator $Y:\mathcal{F}_{n}\to\mathcal{H}$,
and such that $\|[S_{u}y\otimes h_{j}]_{\mathcal{S}}\|<\epsilon$
and $\|[h_{j}\otimes S_{u}y]_{\mathcal{S}}\|<\epsilon$ for any word
$u\in\mathbb{F}_{n}^{+}\mbox{ and }1\leq j\leq q.$
\end{lem}
\begin{proof}
For each $j$, there exists $\overline{h}^{(j)}=(h_{1}^{(j)},...,h_{M}^{(j)})$
in $\mathcal{W}$ such that $h_{j}=h_{1}^{(j)}$. Choose $\epsilon_{0}\in(0,1)$
such that $\epsilon_{0}/(1-\epsilon_{0})<\epsilon$ and $\epsilon_{0}/\|h^{(j)}\|<1$
for $1\leq j\leq q$, and choose $r$ sufficiently large that $\overline{h}^{(j)}$
is orthogonal to $\mathcal{W}_{r}$ for $1\leq j\leq q$.

By Lemma \ref{lem:key-lemma-2}, there exists a unit vector $\overline{x}=(x_{1},...,x_{M})$
in $\mathcal{M}_{r}$ such that $x_{1}=X\xi_{\varnothing}$ for some
intertwining operator $X:\mathcal{F}_{n}\to\mathcal{H}$, and such
that\[
\|x_{1}\|>\max\left\{ 1-\epsilon_{0},\left(1-\frac{\epsilon_{0}^{2}}{\|\overline{h}^{(j)}\|^{2}}\right)^{1/2}:1\leq j\leq q\right\} .\]
 This gives $1/\|x_{1}\|<1/(1-\epsilon_{0})$ and\[
\sum_{i=2}^{M}\|x_{i}\|^{2}=1-\|x_{1}\|^{2}<\frac{\epsilon_{0}^{2}}{\|\overline{h}^{(j)}\|^{2}},\quad1\leq j\leq q.\]
For any word $u$ in $\mathbb{F}_{n}^{+}$, \begin{eqnarray*}
\|[S_{u}x_{1}\otimes h_{1}^{(j)}]_{\mathcal{S}}\| & \leq & \left\Vert \sum_{i=1}^{M}[S_{u}x_{i}\otimes h_{i}^{(j)}]_{\mathcal{S}}\right\Vert +\left\Vert \sum_{i=2}^{M}[S_{u}x_{i}\otimes h_{i}^{(j)}]_{\mathcal{S}}\right\Vert \\
 & = & \left\Vert [S_{u}^{(M)}\overline{x}\otimes\overline{h}^{(j)}]_{\mathcal{S}}\right\Vert +\left\Vert \sum_{i=2}^{M}[S_{u}x_{i}\otimes h_{i}^{(j)}]_{\mathcal{S}}\right\Vert \\
 & = & \left\Vert \sum_{i=2}^{M}[S_{u}x_{i}\otimes h_{i}^{(j)}]_{\mathcal{S}}\right\Vert \\
 & \leq & \left(\sum_{i=2}^{M}\|S_{u}x_{i}\|^{2}\right)^{1/2}\left(\sum_{i=2}^{M}\|h_{i}^{(j)}\|^{2}\right)^{1/2}\\
 & = & \left(\sum_{i=2}^{M}\|x_{i}\|^{2}\right)^{1/2}\left(\sum_{i=2}^{M}\|h_{i}^{(j)}\|^{2}\right)^{1/2}\\
 & < & \frac{\epsilon_{0}}{\|\overline{h}^{(j)}\|}\left(\sum_{i=2}^{M}\|h_{i}^{(j)}\|^{2}\right)^{1/2}\\
 & \leq & \epsilon_{0},\end{eqnarray*}
where we have used the fact that $\overline{x}$ and $\overline{h}^{(j)}$
belong to orthogonal $\mathcal{S}^{(M)}$-invariant subspaces, which
implies that $\|[S_{u}^{(M)}\overline{x}\otimes\overline{h}^{(j)}]_{\mathcal{S}}\|=0$.
Multiplying this inequality by $1/\|x_{1}\|=1/\|S_{u}x_{1}\|$ then
gives \[
\|[S_{u}(x_{1}/\|x_{1}\|)\otimes h_{j}]\|<\epsilon_{0}/(1-\epsilon_{0})\]
 for $1\leq j\leq q$. In the same way we get \[
\|h_{j}\otimes S_{u}(x_{1}/\|x_{1}\|)]\|<\epsilon_{0}/(1-\epsilon_{0})\]
 for $1\leq j\leq q$. Hence we can take $y=x_{1}/\|x_{1}\|$ and
$Y=X/\|x_{1}\|$.
\end{proof}

\begin{lem}
\label{lem:compile-small-approx}Given $h_{1},...,h_{q}$ in $\mathcal{M}$,
$p\geq1$, and $\epsilon>0$, there exists a unit vector $z$ in $\mathcal{M}$
such that\[
\|[S_{u_{1}}z\otimes S_{u_{2}}z]_{\mathcal{S}}-\phi([\xi_{u_{1}}\otimes\xi_{u_{2}}]_{\mathcal{L}_{n}})\|<\epsilon\]
 for all $u_{1}$ and $u_{2}$ in $\mathbb{F}_{n}^{p}$, and such
that $\|[S_{w}z\otimes h_{j}]_{\mathcal{S}}\|<\epsilon$ and $\|[h_{j}\otimes S_{w}z]_{\mathcal{S}}\|<\epsilon$
for all $w\in\mathbb{F}_{n}^{+}$ and $1\leq j\leq q$.
\end{lem}
\begin{proof}
By Lemma \ref{lem:quasi-orthog}, there exists a unit vector $y$
in $\mathcal{M}$ such that $y=Y\xi_{\varnothing}$ for some intertwining
operator $Y:\mathcal{F}_{n}\to\mathcal{H}$, and such that for any
word $w$ in $\mathbb{F}_{n}^{+}$, $\|[S_{w}y\otimes h_{j}]_{\mathcal{S}}\|<\epsilon$
and $\|[h_{j}\otimes S_{w}y]_{\mathcal{S}}\|<\epsilon$ for $1\leq j\leq q$.
By Lemma \ref{lem:word-close-to-wand}, there exists a word $v$ in
$\mathbb{F}_{n}^{+}$ such that $\|[S_{u_{1}}S_{v}y\otimes S_{u_{2}}S_{v}y]_{\mathcal{S}}-\phi([\xi_{u_{1}}\otimes\xi_{u_{2}}]_{\mathcal{L}_{n}})\|<\epsilon$
for any words $u_{1}$ and $u_{2}$ in $\mathbb{F}_{n}^{p}$. Then
$\|[S_{w}S_{v}y\otimes h_{j}]_{\mathcal{S}}\|=\|[S_{wv}y\otimes h_{j}]_{\mathcal{S}}\|<\epsilon$
and $\|[h_{j}\otimes S_{w}S_{v}y]_{\mathcal{S}}\|=\|[h_{j}\otimes S_{wv}y]_{\mathcal{S}}\|<\epsilon$,
so we can take $z=S_{v}y$.
\end{proof}

\begin{lem}
\label{lem:big-approx-1}Given a weak{*}-continuous linear functional
$\pi$ on $\mathcal{S}$, $h_{1},...,h_{q}$ in $\mathcal{M}$, and
$\epsilon>0$, there are vectors $x$ and $y$ in $\mathcal{M}$ such
that 
\begin{enumerate}
\item $\|\pi-[x\otimes y]_{\mathcal{S}}\|<\epsilon$,
\item $\|x\|<(1+\epsilon)\|\pi\|^{1/2}$ and $\|y\|<(1+\epsilon)\|\pi\|^{1/2}$,
\item $\|[x\otimes h_{j}]_{\mathcal{S}}\|<\epsilon$ and $\|[h_{j}\otimes y]_{\mathcal{S}}\|<\epsilon$
for $1\leq j\leq q$.
\end{enumerate}
\end{lem}
\begin{proof}
By scaling $\pi$ and $\epsilon$ if necessary, we can assume that
$\|\pi\|=1$. Choose $\epsilon_{0}>0$ such that $2\epsilon_{0}+3\epsilon_{0}^{2}<\epsilon/2$
and $4\epsilon_{0}+4\epsilon_{0}^{2}<\epsilon+\epsilon^{2}/2$. Since
$\mathcal{L}_{n}$ has property $\mathbb{A}_{1}(1)$, there are vectors
$\xi$ and $\eta$ in $\mathcal{F}_{n}$ such that $[\xi\otimes\eta]_{\mathcal{L}_{n}}=\phi^{-1}(\pi)$,
with $\|\xi\|<1+\epsilon_{0}$ and $\|\eta\|<1+\epsilon_{0}$.

Since $\xi_{\varnothing}$ is cyclic for $\mathcal{L}_{n}$, there
is $p\geq$1 and $C$ and $D$ in the span of $\{L_{u}:u\in\mathbb{F}_{n}^{p}\}$
such that $\|C\xi_{\varnothing}-\xi\|<\epsilon_{0}$ and $\|D\xi_{\varnothing}-\eta\|<\epsilon_{0}$.
Then\[
\|C\xi_{\varnothing}\|\leq\|C\xi_{\varnothing}-\xi\|+\|\xi\|<1+2\epsilon_{0},\]
so $\|C\xi_{\varnothing}\|^{2}<1+\epsilon+\epsilon^{2}/2$, and similarly,
$\|D\xi_{\varnothing}\|^{2}<1+\epsilon+\epsilon^{2}/2$. Also,\begin{eqnarray*}
\|\pi-[C\xi_{\varnothing}\otimes D\xi_{\varnothing}]_{\mathcal{L}_{n}}\| & = & \|[\xi\otimes\eta]_{\mathcal{L}_{n}}-[C\xi_{\varnothing}\otimes D\xi_{\varnothing}]_{\mathcal{L}_{n}}\|\\
 & \leq & \|[(\xi-C\xi_{\varnothing})\otimes\eta]_{\mathcal{L}_{n}}\|+\|[(C\xi_{\varnothing}-\xi)\otimes(\eta-D\xi_{\varnothing})]_{\mathcal{L}_{n}}\|+\\
 &  & \|[\xi\otimes(\eta-D\xi_{\varnothing})]_{\mathcal{L}_{n}}\|\\
 & \leq & \|\xi-C\xi_{\varnothing}\|\|\eta\|+\|\xi-C\xi_{\varnothing}\|\|\eta-D\xi_{\varnothing}\|+\\
 &  & \|\xi\|\|\eta-D\xi_{\varnothing}\|\\
 & < & 2\epsilon_{0}+3\epsilon_{0}^{2}\\
 & < & \epsilon/2.\end{eqnarray*}

Set $A=\Phi(C)$ and $B=\Phi(D)$. If we expand $C$ and $D$ as \[
C=\sum_{u\in\mathbb{F}_{n}^{p}}c_{u}L_{u}\quad\mbox{ and }\quad D=\sum_{u\in\mathbb{F}_{n}^{p}}d_{u}L_{u},\]
then\[
A=\sum_{u\in\mathbb{F}_{n}^{p}}c_{u}S_{u}\quad\mbox{ and }\quad B=\sum_{u\in\mathbb{F}_{n}^{p}}d_{u}S_{u}.\]
Choose $\epsilon_{1}>0$ such that\[
\epsilon_{1}\sum_{u\in\mathbb{F}_{n}^{p}}\left|c_{u}\right|<\epsilon,\quad\epsilon_{1}\sum_{u\in\mathbb{F}_{n}^{p}}\left|\overline{d_{u}}\right|<\epsilon,\quad\epsilon_{1}\sum_{u\in\mathbb{F}_{n}^{p}}\sum_{v\in\mathbb{F}_{n}^{p}}\left|c_{u}\overline{d_{v}}\right|<\epsilon/2,\]
\[
\epsilon_{1}\sum_{u\in\mathbb{F}_{n}^{N}}\sum_{v\in\mathbb{F}_{n}^{N}}\left|c_{u}\overline{c_{v}}\right|<\epsilon+\epsilon^{2}/2,\quad\epsilon_{1}\sum_{u\in\mathbb{F}_{n}^{N}}\sum_{v\in\mathbb{F}_{n}^{N}}\left|d_{u}\overline{d_{v}}\right|<\epsilon+\epsilon^{2}/2.\]
By Lemma \ref{lem:compile-small-approx}, there exists a unit vector
$z$ in $\mathcal{M}$ such that\[
\|[S_{u}z\otimes S_{v}z]_{\mathcal{S}}-\phi([\xi_{u}\otimes\xi_{v}]_{\mathcal{L}_{n}})\|<\epsilon_{1}\]
for any words $u$ and $v$ in $\mathbb{F}_{n}^{p}$, and such that
$\|[S_{u}z\otimes h_{j}]_{\mathcal{S}}\|<\epsilon_{1}$ and $\|[h_{j}\otimes S_{u}z]_{\mathcal{S}}\|<\epsilon_{1}$
for any word $u$ in $\mathbb{F}_{n}^{+}$ and $1\leq j\leq q$. Then

\begin{align*}
\|[Az\otimes Bz]_{\mathcal{S}}\,-\, & \phi([C\xi_{\varnothing}\otimes D\xi_{\varnothing}]_{\mathcal{L}_{n}})\|\\
 & =\left\Vert \sum_{u\in\mathbb{F}_{n}^{p}}\sum_{v\in\mathbb{F}_{n}^{p}}c_{u}\overline{d_{v}}([S_{u}z\otimes S_{v}z]_{\mathcal{S}}-\phi([S_{u}\xi_{\varnothing}\otimes S_{v}\xi_{\varnothing}]_{\mathcal{L}_{n}}))\right\Vert \\
 & \leq\sum_{u\in\mathbb{F}_{n}^{p}}\sum_{v\in\mathbb{F}_{n}^{p}}\left|c_{u}\overline{d_{v}}\right|\|[S_{u}z\otimes S_{v}z]_{\mathcal{S}}-\phi([S_{u}\xi_{\varnothing}\otimes S_{v}\xi_{\varnothing}]_{\mathcal{L}_{n}}))\|\\
 & <\epsilon_{1}\sum_{u\in\mathbb{F}_{n}^{p}}\sum_{v\in\mathbb{F}_{n}^{p}}\left|c_{u}\overline{d_{v}}\right|\\
 & <\epsilon/2.\end{align*}
Hence from above,\[
\|\pi-[Az\otimes Bz]_{\mathcal{S}}\|\leq\|\pi-[\phi([C\xi_{\varnothing}\otimes D\xi_{\varnothing}]_{\mathcal{L}_{n}})\|+\|\phi([C\xi_{\varnothing}\otimes D\xi_{\varnothing}]_{\mathcal{L}_{n}})-[Az\otimes Bz]_{\mathcal{S}}\|<\epsilon.\]

By a similar estimation, \[
\|[Az\otimes Az]_{\mathcal{S}}-\phi([C\xi_{\varnothing}\otimes C\xi_{\varnothing}]_{\mathcal{L}_{n}})\|<\epsilon_{1}\sum_{u\in\mathbb{F}_{n}^{p}}\sum_{v\in\mathbb{F}_{n}^{p}}\left|c_{u}\overline{c_{v}}\right|<\epsilon+\epsilon^{2}/2.\]
Evaluation of these functionals at the identity then implies \[
\epsilon+\epsilon^{2}/2>\|Az\|^{2}-\|C\xi_{\varnothing}\|^{2}\geq\|Az\|^{2}-(1+\epsilon+\epsilon^{2}/2),\]
and hence that $\|Az\|<1+\epsilon.$ In the same way we get $\|Bz\|<1+\epsilon$. 

Finally, \begin{eqnarray*}
\|[Az\otimes h_{j}]_{\mathcal{S}}\| & = & \|\sum_{u\in\mathbb{F}_{n}^{p}}c_{u}[S_{u}z\otimes h_{j}]_{\mathcal{S}}\|\\
 & \leq & \sum_{u\in\mathbb{F}_{n}^{p}}\left|c_{u}\right|\|[S_{u}z\otimes h_{j}]_{\mathcal{S}}\|\\
 & < & \epsilon_{1}\sum_{u\in\mathbb{F}_{n}^{p}}\left|c_{u}\right|\\
 & < & \epsilon,\end{eqnarray*}
and in the same way we get \[
\|[h_{j}\otimes Bz]_{\mathcal{S}}\|<\epsilon_{1}\sum_{u\in\mathbb{F}_{n}^{p}}\left|\overline{d}_{u}\right|<\epsilon.\]
Hence we can take $x=Az$ and $y=Bz$.
\end{proof}

\begin{thm}
Given a weak{*}-continuous linear functional $\pi$ on $\mathcal{S}$
and $\epsilon>0$, there are vectors $x$ and $y$ in $\mathcal{H}$
such that $\pi=[x\otimes y]_{\mathcal{S}}$, $\|x\|<(1+\epsilon)\|\pi\|^{1/2}$,
and $\|y\|<(1+\epsilon)\|\pi\|^{1/2}$. In other words, $\mathcal{S}$
has property $\mathbb{A}_{1}(1)$.
\end{thm}
\begin{proof}
By scaling $\pi$ if necessary, we can assume that $\|\pi\|=1$. Choose
$\alpha>0$ such that $(1+\alpha)/(1-\alpha)<1+\epsilon$. Note that
$\alpha^{k}\to0$ as $k\to\infty$. We claim that for $k\geq1$, we
can find $x_{k}$ and $y_{k}$ in $\mathcal{M}$ such that
\begin{enumerate}
\item $\|\pi-[x_{k}\otimes y_{k}]_{\mathcal{S}}\|<\alpha^{2k}$,
\item $\|x_{k}\|<(1+\alpha)(1+\alpha+...+\alpha^{k-1})$ and $\|y_{k}\|<(1+\alpha)(1+\alpha+...+\alpha^{k-1})$,
\item $\|x_{k}-x_{k-1}\|<(1+\alpha)\alpha^{k-1}$ and $\|y_{k}-y_{k-1}\|<(1+\alpha)\alpha^{k-1}$
for $k\geq2$.
\end{enumerate}
Setting $x_{0}=0$ and $y_{0}=0$, Lemma \ref{lem:big-approx-1} easily
implies this is true for $k=1$. Proceeding by induction, suppose
that we have found $x_{k}$ and $y_{k}$ satisfying these conditions.
Choose $\epsilon_{0}>0$ such that $\epsilon_{0}<\alpha$ and $\epsilon_{0}<\alpha^{2(k+1)}/3$.
By Lemma \ref{lem:big-approx-1}, there are $x'$ and $y'$ in $\mathcal{M}$
such that
\begin{enumerate}
\item $\|\pi-[x_{k}\otimes y_{k}]_{\mathcal{S}}-[x'\otimes y']_{\mathcal{S}}\|<\epsilon_{0}$,
\item $\|x'\|<(1+\epsilon_{0})\|\pi-[x_{k}\otimes y_{k}]_{\mathcal{S}}\|^{1/2}$
and\\
 $\|y'\|<(1+\epsilon_{0})\|\pi-[x_{k}\otimes y_{k}]_{\mathcal{S}}\|^{1/2}$,
\item $\|[x'\otimes y_{k}]_{\mathcal{S}}\|<\epsilon_{0}$ and $\|[x_{k}\otimes y']_{\mathcal{S}}\|<\epsilon_{0}$.
\end{enumerate}
Set $x_{k+1}=x_{k}+x'$, and $y_{k+1}=y_{k}+y'$. Then \begin{eqnarray*}
\|\pi-[x_{k+1}\otimes y_{k+1}]_{\mathcal{S}}\| & = & \|\pi-[(x_{k}+x')\otimes(y_{k}+y')]_{\mathcal{S}}\|\\
 & \leq & \|\pi-[x_{k}\otimes y_{k}]_{\mathcal{S}}-[x'\otimes y']_{\mathcal{S}}\|+\|[x_{k}\otimes y']_{\mathcal{S}}\|+\\
 &  & \|[x'\otimes y_{k}]_{\mathcal{S}}\|\\
 & < & 3\epsilon_{0}\\
 & < & \alpha^{2(k+1)}.\end{eqnarray*}
Also, \[
\|x'\|<(1+\epsilon_{0})\|\pi-[x_{k}\otimes y_{k}]_{\mathcal{S}}\|^{1/2}<(1+\alpha)\alpha^{k},\]
which gives\[
\|x_{k+1}\|=\|x_{k}+x'\|\leq\|x_{k}\|+\|x'\|<(1+\alpha)(1+\alpha+...+\alpha^{k})\]
and \[
\|x_{k+1}-x_{k}\|=\|x'\|<(1+\alpha)\alpha^{k}.\]
Symmetrically, $\|y_{k+1}\|<(1+\alpha)(1+\alpha+...+\alpha^{k})$
and $\|y_{k+1}-y_{k}\|<(1+\alpha)\alpha^{k}$, which establishes the
claim.

Now for $l>k$, \begin{eqnarray*}
\|x_{l}-x_{k}\| & \leq & \|x_{l}-x_{l-1}\|+...+\|x_{k+1}-x_{k}\|\\
 & < & (1+\alpha)(\alpha^{l-1}+...+\alpha^{k})\\
 & \leq & \alpha^{l-1}(1+\alpha)/(1-\alpha),\end{eqnarray*}
so the sequence $(x_{k})$ is Cauchy. Let $x=\lim_{k}x_{k}$. Then\[
\|x\|=\lim_{k}\|x_{k}\|\leq\lim_{k}(1+\alpha)(1+\alpha+...+\alpha^{k-1})=(1+\alpha)/(1-\alpha)<1+\epsilon.\]
Similarly, the sequence $\{y_{k}\}$ is Cauchy. Letting $y=\lim_{k}y_{k}$
be its limit, $\|y\|<1+\epsilon$. Finally, we have\[
\|\pi-[x\otimes y]_{\mathcal{S}}\|=\lim\|\pi-[x_{k}\otimes y_{k}]_{\mathcal{S}}\|\leq\lim\alpha^{2k}=0,\]
so $\pi=[x\otimes y]_{\mathcal{S}}$.
\end{proof}

\begin{thm}
\label{thm:wand-vec}Every type L free semigroup algebra has a wandering
vector.
\end{thm}
\begin{proof}
Let $\mathcal{S}$ be a type L free semigroup algebra, and let $\mathcal{S}_{0}$
denote the weak-operator-closed ideal generated by $S_{1},...,S_{n}$.
Since $\mathcal{S}$ is type L, $\mathcal{S}_{0}$ is proper, and
in particular doesn't contain the identity. Let $\pi_{0}$ denote
the weak-operator continuous linear functional which annihilates $\mathcal{S}_{0}$
and satisfies $\pi(I)=1$.

Since $\mathcal{S}$ has property $\mathbb{A}_{1}(1)$, there are
vectors $x$ and $y$ in $\mathcal{H}$ such that $\pi_{0}(A)=(Ax,y)$
for all $A$ in $\mathcal{S}$. This implies $(S_{w}x,y)=0$ for all
$w\in\mathbb{F}_{n}^{+}\backslash\{\varnothing\}$, so $y$ is orthogonal
to the subspace $\mathcal{S}_{0}[x]$. However, $(x,y)=\pi(I)=1$,
so $y$ is not orthogonal to the subspace $\mathcal{S}[x]$. Hence
$\mathcal{S}[x]\ominus\mathcal{S}_{0}[x]$ is nonempty.

Let $z$ be a unit vector in $\mathcal{S}[x]\ominus\mathcal{S}_{0}[x]$.
Then the subspace $\mathcal{S}_{0}[z]$ is contained in the subspace
$\mathcal{S}_{0}[x]$, and in particular, is orthogonal to $z$. Hence
$(S_{w}z,z)=0$ for all $w\in\mathbb{F}_{n}^{+}\backslash\{\varnothing\}$.
Let $u$ and $v$ be distinct words in $\mathbb{F}_{n}^{+}$ such
that $|u|\leq|v|$. Then $S_{u}^{*}S_{v}$ is in $\mathcal{S}_{0}$,
so $(S_{u}z,S_{v}z)=(z,S_{u}^{*}S_{v}z)=0$. By symmetry, it follows
that $(S_{u}z,S_{v}z)=0$ for every pair of distinct words $u$ and
$v$ in $\mathbb{F}_{n}^{+}$. Thus $z$ is a wandering vector for
$\mathcal{S}$. 
\end{proof}

\begin{cor}
A free semigroup algebra is either a von Neumann algebra, or it contains
a wandering vector.
\end{cor}
\begin{proof}
Let $\mathcal{S}$ be a free semigroup algebra. By the general structure
theorem for free semigroup algebras \cite{MR1823866}, $\mathcal{S}$
is either a von Neumann algebra, or it has a type L part. In the latter
case, by Theorem \ref{thm:wand-vec}, $\mathcal{S}$ has a wandering
vector. 
\end{proof}

By Theorem 4.1 of \cite{MR2139108}, every free semigroup algebra
which has a wandering vector is reflexive. Thus we have established
the following result. 

\begin{cor}
\label{cor:reflexive}Every free semigroup algebra is reflexive.
\end{cor}

Theorem 4.2 of \cite{MR2139108} shows that every type L free semigroup
algebra which has a wandering vector is hyper-reflexive with hyper-reflexivity
constant at most $55$. This gives the following result, which we
will refine in section \ref{sec:hyper-reflexivity}. 

\begin{cor}
\label{cor:type-l-hyper-ref-55}Every type L free semigroup algebra
is hyper-reflexive with hyper-reflexivity constant at most $55$.
\end{cor}

\section{\label{sec:hyper-reflexivity}Hyper-reflexivity and the Factorization
of Linear Functionals}

In this section we will show that the predual of every type L free
semigroup algebra with $n\geq2$ generators satisfies a very strong
factorization property. By a result of Bercovici \cite{MR1641578},
we will obtain as a consequence that every such algebra is hyper-reflexive
with hyper-reflexivity constant at most $3$.

\begin{defn}
A weak{*}-closed subspace $\mathcal{S}$ of $\mathcal{B}(\mathcal{H})$
is said to have property $\mathcal{X}_{0,1}$ if given a weak{*}-continuous
linear functional $\pi$ on $\mathcal{S}$ with $\|\pi\|\leq1$, $h_{1},...,h_{q}$
in $\mathcal{H}$, and $\epsilon>0$, there are vectors $x$ and $y$
in $\mathcal{H}$ such that
\begin{enumerate}
\item $\|\pi-[x\otimes y]_{\mathcal{S}}\|<\epsilon$,
\item $\|x\|\leq1$ and $\|y\|\leq1$,
\item $\|[x\otimes h_{j}]_{\mathcal{S}}\|<\epsilon$ and $\|[h_{j}\otimes y]_{\mathcal{S}}\|<\epsilon$
for $1\leq j\leq q$.
\end{enumerate}
\end{defn}

Bercovici \cite{MR1641578} showed that any weak{*}-closed algebra
whose commutant contains two isometries with pairwise orthogonal ranges
has property $\mathcal{X}_{0,1}$, and showed that any weak{*}-closed
algebra with propety $\mathcal{X}_{0,1}$ is hyper-reflexive with
hyper-reflexivity constant at most $3$. For $n\geq2$, this includes
$\mathcal{L}_{n}$. We will show that every type L free semigroup
algebra with $n\geq2$ generators has property $\mathcal{X}_{0,1}$.

We require the following result which is implied by Lemma 1.2 in \cite{MR1823067}.

\begin{lem}
\label{lem:strongly}Given an isometry $V$ in $\mathcal{R}_{n}$,
vectors $\nu_{1},...,\nu_{q}$ in $\mathcal{F}_{n}$, and $\epsilon>0$,
there exists $m$ such that $\|(V^{*})^{m}\eta_{j}\|<\epsilon$ for
$1\leq j\leq q$.
\end{lem}

For the remainder of this section we fix a type L free semigroup algebra
$\mathcal{S}$ with $n\geq2$ generators acting on a Hilbert space
$\mathcal{H}$, and we let $\mathcal{Z}$ denote the weak{*} closure
of $\mathcal{S}+\mathcal{S}^{*}$. Let $\Phi$ denote the canonical
map from $\mathcal{L}_{n}$ to $\mathcal{S}$. By Theorem \ref{thm:exten-canon},
we can extend $\Phi$ to a map from the set $\mathcal{T}_{R}$ of
$R$-Toeplitz operators to $\mathcal{Z}$, and this extension is a
complete isometry and a weak{*}-to-weak{*} homeomorphism. 

For $x$ and $y$ in $\mathcal{H}$, we will need to take care to
distinguish between the weak-operator-continuous vector functional
$[x\otimes y]_{\mathcal{S}}$ defined on $\mathcal{S}$, and the weak-operator-continuous
vector functional $[x\otimes y]_{\mathcal{Z}}$ defined on $\mathcal{Z}$. 

The following lemma is a variation of an argument of Bercovici \cite{MR1641578}.
It was kindly provided by Ken Davidson.

\begin{lem}
\label{lem:bercovici}Given isometries $U$ and $V$ in $\mathcal{R}_{n}$
with orthogonal ranges, vectors $\xi$ and $\nu$ in $\mathcal{F}_{n}$
with $\nu$ in the kernel of $U^{*}$, and $\epsilon>0$, define\[
\eta_{k}=\frac{1}{\sqrt{k}}\sum_{i=1}^{k}U^{i}V\xi.\]
Then $\lim_{k}\|[\nu\otimes\eta_{k}]_{\mathcal{T}_{R}}\|=0.$ 
\end{lem}
\begin{proof}
Let $H^{2}$ denote the Hardy-Hilbert space with orthonormal basis
$\{e_{k}:k\geq0$\}. For $k\geq0$, define $Y:H^{2}\to\mathcal{F}_{n}$
by $Ye_{k}=U^{k}V\xi$ and $Z:H^{2}\to\mathcal{F}_{n}$ by $Ze_{k}=U^{k}\nu$
for $k\geq0$. Note that $Y$ and $Z$ are isometries. For $T$ in
$\mathcal{T}_{R}$, by Lemma \ref{lem:factor-r-toeplitz} we can factor
$T$ as $T=A^{*}B$, for $A$ and $B$ in $\mathcal{L}_{n}$. Then
\begin{eqnarray*}
(Y^{*}TZe_{j},e_{i}) & = & (A^{*}BU^{j}\nu,U^{i}V\xi)\\
 & = & (A^{*}V^{*}(U^{*})^{i}U^{j}B\nu,\xi)\\
 & = & \begin{cases}
0 & \mbox{if }i<j\\
c_{i-j} & \mbox{if }i\geq j\end{cases},\end{eqnarray*}
where $c_{i-j}=(A^{*}B\nu,U^{i-j}V\xi)=(T\nu,U^{i-j}V\xi)$. This
implies that $Y^{*}TZ$ is an analytic Toeplitz operator with symbol
$f$, for some $f$ in $H^{\infty}$. Note that $\|f\|_{\infty}=\|Y^{*}TZ\|\leq\|T\|$.
Hence\begin{eqnarray*}
\left|(T,[\nu\otimes\eta_{k}]_{\mathcal{T}_{R}})\right| & = & \left|(T\nu,\frac{1}{\sqrt{k}}\sum_{i=1}^{k}U^{i}V\xi)\right|\\
 & = & \frac{1}{\sqrt{k}}\left|\sum_{i=1}^{k}(T\nu,U^{i}V\xi)\right|\\
 & = & \frac{1}{\sqrt{k}}\left|\sum_{i=1}^{k}c_{i}\right|\\
 & \leq & \frac{1}{\sqrt{k}}\|D_{k}\|_{1}\|f\|_{\infty}\\
 & \leq & \frac{1}{\sqrt{k}}\|D_{k}\|_{1}\|T\|,\end{eqnarray*}
where $\|D_{k}\|_{1}$ denotes the $L^{1}$-norm of the Dirichlet
kernel. Using the well-known fact that $\|D_{k}\|_{1}$ grows logarithmically
as $k\to\infty$ gives $\lim_{k}\|[\nu\otimes\eta_{k}]_{\mathcal{T}_{R}}\|=0$.
\end{proof}

\begin{lem}
\label{lem:big-approx}Given vectors $h_{1},...,h_{q}$ in $\mathcal{H}$
and $\epsilon>0$, there exists an intertwining operator $Y:\mathcal{F}_{n}\to\mathcal{H}$
such that $\|Y\xi_{\varnothing}\|=1$ and $\|[Y\xi_{\varnothing}\otimes h_{i}]_{\mathcal{Z}}\|<\epsilon$\textup{
for $1\leq i\leq q$.}
\end{lem}
\begin{proof}
For $1\leq i\leq q$, let $H_{i}:\mathcal{F}_{n}\to\mathcal{H}$ be
an intertwining operator such that $\|H_{i}\xi_{\varnothing}-h_{i}\|<\epsilon/2$.
Since $\mathcal{S}$ is type L, by Theorem \ref{thm:wand-vec} there
is an isometric intertwining operator $X:\mathcal{F}_{n}\to\mathcal{H}$.
Then each $H_{i}^{*}X$ is an L-Toeplitz operator, so by Lemma \ref{lem:factor-r-toeplitz},
we can write $H_{i}^{*}X=A_{i}^{*}B_{i}$, for some $A_{i}$ and $B_{i}$
in $\mathcal{R}_{n}$ such that $A_{i}$ and $B_{i}$ are bounded
below. Let $C_{i}=R_{1^{i}2}B_{i}$, and let $D=\sum_{i=1}^{k}R_{1^{i}2}A_{i}$.
Then $D$ is bounded below and $H_{i}^{*}X=C_{i}^{*}D$. Using inner-outer
factorization, write $D=UF$ for $U$ and $F$ in $\mathcal{R}_{n}$,
where $U$ is inner and $F$ is outer. Then $F$ is bounded below
since $D$ is, and hence is invertible.

By Lemma \ref{lem:strongly}, there exists $m$ such that $\|(U^{*})^{m}C_{i}\xi_{\varnothing}\|<\epsilon/(8\|F\|)$
for $1\leq i\leq q$. Write $C_{i}\xi_{\varnothing}=\nu_{i}+\omega_{i}$,
where $\|\omega_{i}\|<\epsilon/(8\|F\|)$, and $\nu_{i}$ is in the
kernel of $(U^{*})^{m}$. Set $V=U^{m}R_{1}$ and $W=U^{m}R_{2}$.
Then $V$ and $W$ are isometries in $\mathcal{R}_{n}$ with pairwise
orthogonal ranges. Note that $\nu_{i}$ is in the kernel of $V^{*}$.
For $k\geq1$, define intertwining operators $Y_{k}:\mathcal{F}_{n}\to\mathcal{H}$
by \[
Y_{k}=XF^{-1}\frac{1}{\sqrt{k+1}}\sum_{j=0}^{k}U^{m-1}R_{1}V^{j}W,\]
and define \[
\eta_{k}=\frac{1}{\sqrt{k}}\sum_{j=1}^{k}V^{j}W\xi_{\varnothing}.\]
Note that $\eta_{k}$ is a unit vector.

Using the fact that $V=DF^{-1}U^{m-1}R_{1}$, we compute \begin{eqnarray*}
H_{i}^{*}Y_{k} & = & H_{i}^{*}XF^{-1}\frac{1}{\sqrt{k+1}}\sum_{j=0}^{k}U^{m-1}R_{1}V^{j}W\\
 & = & C_{i}^{*}DF^{-1}\frac{1}{\sqrt{k+1}}\sum_{j=0}^{k}U^{m-1}R_{1}V^{j}W\\
 & = & C_{i}^{*}\frac{1}{\sqrt{k+1}}\sum_{j=1}^{k+1}V^{j}W.\end{eqnarray*}
Then for $T$ in $\mathcal{T}_{R}$, \begin{eqnarray*}
(TH_{i}\xi_{\varnothing},Y_{k}\xi_{\varnothing}) & = & (\Phi^{-1}(T)C_{i}\xi_{\varnothing},\frac{1}{\sqrt{k+1}}\sum_{j=1}^{k+1}V^{j}W\xi_{\varnothing})\\
 & = & (\Phi^{-1}(T)C_{i}\xi_{\varnothing},\eta_{k+1}),\end{eqnarray*}
Hence $\|[H_{i}\xi_{\varnothing}\otimes Y_{k}\xi_{\varnothing}]_{\mathcal{Z}}\|=\|[C_{i}\xi_{\varnothing}\otimes\eta_{k+1}]_{\mathcal{T}_{R}}\|.$
By Lemma \ref{lem:bercovici}, we can choose $r$ sufficiently large
that $\|[\nu_{i}\otimes\eta_{r+1}]_{\mathcal{T}_{R}}\|<\epsilon/(8\|F\|)$.
This gives\begin{eqnarray*}
\|[C_{i}\xi_{\varnothing}\otimes\eta_{r+1}]_{\mathcal{T}_{R}}\| & \leq & \|[\nu_{i}\otimes\eta_{r+1}]_{\mathcal{T}_{R}}\|+\|[\omega_{i}\otimes\eta_{r+1}]_{\mathcal{T}_{R}}\|\\
 & \leq & \|[\nu_{i}\otimes\eta_{r+1}]_{\mathcal{T}_{R}}]\|+\|\omega_{i}\|\|\eta_{r+1}\|\\
 & < & \epsilon/(4\|F\|).\end{eqnarray*}
Thus $\|[H_{i}\xi_{\varnothing}\otimes Y_{k}\xi_{\varnothing}]_{\mathcal{Z}}\|<\epsilon/(4\|F\|)$
for $1\leq i\leq q$.

Now, \begin{eqnarray*}
\|Y_{r}\xi_{\varnothing}\|^{2} & = & \|XF^{-1}\frac{1}{\sqrt{r+1}}\sum_{j=0}^{r}U^{m-1}R_{1}V^{j}W\xi_{\varnothing}\|^{2}\\
 & \geq & \frac{1}{(r+1)\|F\|^{2}}\|\sum_{j=0}^{r}U^{m-1}R_{1}V^{j}W\xi_{\varnothing}\|^{2}\\
 & = & \frac{r}{(r+1)\|F\|^{2}},\end{eqnarray*}
which implies \[
\|Y_{r}\xi_{\varnothing}\|\geq\frac{1}{2\|F\|}.\]
Setting $Y=Y_{r}/\|Y_{r}\xi_{\varnothing}\|$, it follows that\begin{eqnarray*}
\|[H_{i}\xi_{\varnothing}\otimes Y\xi_{\varnothing}]_{\mathcal{Z}}\| & = & \frac{1}{\|Y_{p}\xi_{\varnothing}\|}\|[H_{i}\xi_{\varnothing}\otimes Y_{p}\xi_{\varnothing}]_{Z}\|\\
 & \leq & 2\|F\|\|[H_{i}\xi_{\varnothing}\otimes Y_{k}\xi_{\varnothing}]_{\mathcal{Z}}\|\\
 & < & \epsilon/2.\end{eqnarray*}
Thus\begin{eqnarray*}
\|[h_{j}\xi_{\varnothing}\otimes Y\xi_{\varnothing}]_{Z}\| & \leq & \|[(h_{j}-H_{j}\xi_{\varnothing})\otimes Y\xi_{\varnothing}]_{\mathcal{Z}}\|+\|[H_{j}\xi_{\varnothing}\otimes Y\xi_{\varnothing}]_{\mathcal{Z}}\|\\
 & \leq & \|h_{j}-H_{j}\xi_{\varnothing}\|\|Y\xi_{\varnothing}\|+\|[H_{j}\xi_{\varnothing}\otimes Y\xi_{\varnothing}]_{\mathcal{Z}}\|\\
 & < & \epsilon.\end{eqnarray*}

\end{proof}

\begin{lem}
\label{lem:x01-key-lemma}Given vectors $h_{1},...,h_{q}$ in $\mathcal{H}$,
$p\geq1$, and $\epsilon>0$, there exists a unit vector $z$ in $\mathcal{H}$
such that \[
\|[S_{u_{1}}z\otimes S_{u_{2}}z]_{\mathcal{S}}-\phi([\xi_{u_{1}}\otimes\xi_{u_{2}}]_{\mathcal{L}_{n}})\|<\epsilon\]
 for all $u_{1}$ and $u_{2}$ in $\mathbb{F}_{n}^{p}$, and such
that $\|[S_{w}z\otimes h_{i}]_{\mathcal{S}}\|<\epsilon$ and $\|[h_{i}\otimes S_{w}z]_{\mathcal{S}}\|<\epsilon$
for all $w\in\mathbb{F}_{n}^{+}$ and $1\leq i\leq q$.
\end{lem}
\begin{proof}
By Lemma \ref{lem:big-approx}, there is an intertwining operator
$Y:\mathcal{F}_{n}\to\mathcal{H}$ such that $\|Y\xi_{\varnothing}\|=1$
and $\|[Y\xi_{\varnothing}\otimes h_{i}]_{\mathcal{Z}}\|<\epsilon$
for $1\leq i\leq q$. By Lemma \ref{lem:word-close-to-wand}, there
is a word $v$ in $\mathbb{F}_{n}^{+}$ such that $\|[S_{u_{1}v}Y\xi_{\varnothing}\otimes S_{u_{2}v}Y\xi_{\varnothing}]_{\mathcal{S}}-\phi([\xi_{u_{1}}\otimes\xi_{u_{2}}]_{\mathcal{L}_{n}})\|<\epsilon$
for all words $u_{1}$and $u_{2}$ in $\mathbb{F}_{n}^{p}$. Set $z=S_{v}Y\xi_{\varnothing}$.

For $T$ in $\mathcal{Z}$ and $w\in\mathbb{F}_{n}^{+}$,\begin{eqnarray*}
|(T,[S_{w}z\otimes h_{j}]_{Z})| & = & |(TS_{w},[z\otimes h_{j}]_{\mathcal{Z}})|\\
 & \leq & \|TS_{w}\|\|[z\otimes h_{j}]_{\mathcal{Z}}\|\\
 & \leq & \|T\|\|[z\otimes h_{J}]_{\mathcal{Z}}\|.\end{eqnarray*}
Hence $\|[S_{w}z\otimes h_{i}]_{\mathcal{Z}}\|\leq\|[z\otimes h_{i}]_{\mathcal{Z}}\|<\epsilon$
and similarly, $\|[h_{i}\otimes S_{w}z]_{Z}\|\leq\|[h_{i}\otimes z]_{\mathcal{Z}}\|<\epsilon$.
In particular, restricting to $\mathcal{S}$ gives $\|[S_{w}z\otimes h_{i}]_{\mathcal{S}}\|<\epsilon$
and $\|[h_{i}\otimes S_{w}z]_{\mathcal{S}}\|<\epsilon$.
\end{proof}

Lemma \ref{lem:x01-key-lemma} is essentially a strengthened version
of Lemma \ref{lem:compile-small-approx}, in the sense that the $h_{i}$'s
in the hypothesis can be completely arbitrary.

\begin{lem}
\label{lem:refined-key-lemma}Given a weak{*}-continuous linear functional
$\pi$ on $\mathcal{S}$, $h_{1},...,h_{q}$ in $\mathcal{H}$, and
$\epsilon>0$, there are vectors $x$ and $y$ in $\mathcal{H}$ such
that 
\begin{enumerate}
\item $\|\pi-[x\otimes y]_{\mathcal{S}}\|<\epsilon$,
\item $\|x\|<(1+\epsilon)\|\pi\|^{1/2}$ and $\|y\|<(1+\epsilon)\|\pi\|^{1/2}$,
\item $\|[x\otimes h_{j}]_{\mathcal{S}}\|<\epsilon$ and $\|[h_{j}\otimes y]_{\mathcal{S}}\|<\epsilon$
for $1\leq j\leq q$.
\end{enumerate}
\end{lem}
\begin{proof}
The proof follows exactly as in the proof of Lemma 4.11, using Lemma
\ref{lem:x01-key-lemma} in place of Lemma 4.10.
\end{proof}

Lemma \ref{lem:refined-key-lemma} clearly implies the desired result.

\begin{thm}
Every type L free semigroup algebra with $n\geq2$ generators has
property $\mathcal{X}_{0,1}$. 
\end{thm}
\begin{cor}
Every type L free semigroup algebra with $n\geq2$ generators is hyper-reflexive
with hyper-reflexivity constant at most $3$.
\end{cor}

\begin{rem}
We note here that the results in Section \ref{sec:hyper-reflexivity}
can be proved independently of Section \ref{sec:wandering}. Let $\mathcal{S}$
be a type L free semigroup algebra. By Lemma \ref{thm:amp-wand-vec},
for some $m$, the ampliation $\mathcal{S}^{(m)}$ has a wandering
vector. By the results in this section, $\mathcal{S}^{(m)}$ has property
$\mathcal{X}_{0,1}$. By Theorem 3.8 of \cite{MR787041}, it follows
that $\mathcal{S}$ has property $\mathcal{X}_{(m-1/m),1}$, which
by \cite{MR1641578} implies that $\mathcal{S}$ is hyper-reflexive,
and in particular is reflexive. It follows that $\mathcal{S}$ has
a wandering vector, and we can apply the results of Section \ref{sec:hyper-reflexivity}
to obtain that $\mathcal{S}$ has property $\mathcal{X}_{0,1}$.
\end{rem}

\section{Concluding Remarks}

In \cite{MR1823866}, Davidson, Katsoulis, and Pitts posed the following
four questions about free semigroup algebras.

\begin{enumerate}
\item Can a free semigroup algebra be a von Neumann algebra? 
\item Does every type L free semigroup algebra have a wandering vector?
\item Is every free semigroup algebra reflexive, or even hyper-reflexive?
\item Is the restriction of a type L free semigroup algebra to an invariant
subspace also of type L?
\end{enumerate}

In \cite{MR2186356} (see also \cite{MR2204288}), Read answered (1)
in the affirmative by showing that $\mathcal{B}(\mathcal{H})$ is
a free semigroup algebra. Theorem \ref{thm:wand-vec} gives an affirmative
answer to (2), and Corollary \ref{cor:reflexive} gives an affirmative
answer to the first part of (3). Corollary \ref{cor:type-l-hyper-ref-55}
partially answers the second part of (3) in the affirmative, as we
now explain.

Recall that by the general structure theorem for free semigroup algebras
\cite{MR1823866}, every free semigroup algebra decomposes into $2\times2$
block-lower-triangular form, where the left column is a slice of a
von Neumann algebra, and the bottom right entry is a type L free semigroup
algebra. By Corollary \ref{cor:type-l-hyper-ref-55}, we know that
the type L part is hyper-reflexive. The difficulty in proving the
hyper-reflexivity of the entire algebra comes down to the fact that
it is an open question whether every von Neumann algebra is hyper-reflexive.

For $n=1$, (4) has a negative answer, as the following example from
\cite{MR1823866} shows.

\begin{example}
\label{exa:abs-con-not-type-l}Let $U$ denote the bilateral shift,
i.e. the operator of multiplication by $z$ on $L^{2}(\mathbb{T},\mu)$,
where $\mu$ is Lebesgue measure. Let $\mathbf{1}_{[0,\pi]}$ and
$\mathbf{1}_{[\pi,2\pi]}$ denote the characteristic function for
the intervals $[0,\pi]$ and $[\pi,2\pi]$ respectively, and define
measures $\mu_{1}$ and $\mu_{2}$ by $\mu_{1}=\mathbf{1}_{[0,\pi]}\mu$
and $\mu_{2}=\mathbf{1}_{[\pi,2\pi]}\mu$. Let $U_{1}$ and $U_{2}$
be the operators of multiplication by $z$ on $L^{2}(\mathbb{T},\mu_{1})$
and $L^{2}(\mathbb{T},\mu_{2})$ respectively. Since the support of
$\mu_{1}$ and $\mu_{2}$ are proper measurable subsets of the circle,
by \cite{MR0048700}, the weak-operator-closed algebras generated
by $U_{1}$ and $U_{2}$ are self-adjoint. On the other hand, $U_{1}\oplus U_{2}$
is unitarily equivalent to the bilateral shift $U$, so the weak-operator-closed
algebra it generates is isomorphic to the analytic Toeplitz algebra.
\end{example}

For $n\geq2$, the answer to (4) is unknown, but it is related to
the notion of absolute continuity, which we now discuss.

The \emph{non-commutative analytic disk algebra} $\mathcal{A}_{n}$
is the norm-closed (non-self-adjoint) algebra generated by $L_{1},...,L_{n}$.
Popescu \cite{MR1343719} showed that the norm-closed algebra generated
by any row isometry of size $n$ is completely isometrically isomorphic
to $\mathcal{A}_{n}$. 

\begin{defn}
\label{def:abs-con}Let $\sigma$ be a $*$-extendible representation
of $\mathcal{A}_{n}$ on a Hilbert space $\mathcal{H}$. These are
precisely the representations with the property that $[\sigma(S_{1})\enspace\cdots\enspace\sigma(S_{n})]$
is a row isometry. We say that $\sigma$ is \emph{absolutely continuous}
if, for every $x$ in $\mathcal{H}$, the linear functional on $\mathcal{A}_{n}$
given by\[
A\to(\sigma(A)x,x),\quad A\in\mathcal{A}_{n}\]
extends to a weak{*}-continuous linear functional on $\mathcal{L}_{n}$.
\end{defn}

If a representation $\sigma$ is type L, i.e. if $\sigma(L_{1}),...,\sigma(L_{n})$
generate a type L free semigroup algebra, then the fact that every
type L free semigroup algebra is completely isometrically isomorphic
and weak{*}-to-weak{*} homeomorphic to $\mathcal{L}_{n}$ immediately
implies that $\sigma$ is absolutely continuous. In \cite{MR2139108},
Davidson, Li, and Pitts conjectured that for $n\geq2$, the converse
of this holds. In other words, they conjectured that for $n\geq2$,
every absolutely continuous representation of $\mathcal{A}_{n}$ is
actually of type L.

For $n=1$, this is not true. Indeed, it's clear from Definition \ref{def:abs-con}
that a representation which is a direct summand of an absolutely continuous
representation is also absolutely continuous. It follows that $U_{1}$
and $U_{2}$ from Example \ref{exa:abs-con-not-type-l} correspond
to absolutely continuous representations. However, the weak-operator-closed
algebras generated by $U_{1}$ and $U_{2}$ are self-adjoint, and
in particular are not of type L.

A proof of the conjecture would therefore provide an interesting example
of a result in the non-commutative setting which has no commutative
counterpart. On the other hand, it would be even more interesting
if this conjecture was false, as we now explain.

It was shown in \cite{MR2383708} that for $n\geq2$, if $\sigma$
is an absolutely continuous representation of $\mathcal{A}_{n}$,
then the infinite ampliation $\sigma^{(\infty)}$ is a representation
of type L. In this case, the weak{*}-closed algebra generated by $\sigma(L_{1}),...,\sigma(L_{n})$
is algebraically isomorphic to $\mathcal{L}_{n}$, and so is non-self-adjoint,
while the weak-operator-closed algebra generated by $\sigma(L_{1}),...,\sigma(L_{n})$
is self-adjoint. This is not completely implausible, as Loebl and
Muhly \cite{MR504460} have demonstrated the existence of weak{*}-closed
non-self-adjoint algebras of analytic type whose weak closure is self-adjoint.

Our results seem to provide some evidence that this conjecture is
true. Davidson, Li, and Pitts \cite{MR2139108} showed that if an
absolutely continuous representation of $\mathcal{A}_{n}$ has a wandering
vector, then it is type L. On the other hand, we know by Theorem \ref{thm:wand-vec}
that if a representation of $\mathcal{A}_{n}$ is type L, then it
has a wandering vector. Therefore, the truth of the conjecture is
equivalent to the existence of a wandering vector for every absolutely
continuous representation. 

For $n\geq2$, suppose that $\sigma$ is an absolutely continuous
representation of $\mathcal{A}_{n}$. One would like to use the methods
in the present paper to find a wandering vector for $\sigma$. Unfortunately,
we no longer know that some finite ampliation of $\sigma$ is type
L, which means we can't make use of Theorem \ref{thm:amp-wand-vec}.
However, most of our approximation techniques do still work in this
setting, so it seems plausible that some modification of our methods
which avoids the dependence on Theorem \ref{thm:amp-wand-vec} could
yield a wandering vector for $\sigma$.

\begin{acknowledgement*}
The author is grateful to Ken Davidson for his advice and support.
\end{acknowledgement*}
\bibliographystyle{amsplain}
\bibliography{wandering}

\end{document}